\documentclass[11pt,lettersize,leqno,fullpage]{article}
\usepackage[usenames,dvipsnames]{color}
\usepackage{amsmath, amssymb, amscd, amsthm, amsfonts, amsbsy, appendix}
\usepackage{bbm, stmaryrd, mathrsfs, bm, linearb, upgreek, skak, textcomp, mathdots}
\usepackage{mathtools, titlesec, extarrows, tocloft, url, dsfont}
\usepackage[intoc]{nomencl}
\usepackage[nottoc]{tocbibind}
\usepackage[all, cmtip]{xy} 
\usepackage{graphicx, pgf}
\definecolor{AliceBlue}{rgb}{0.94,0.97,1.00}
\definecolor{AntiqueWhite1}{rgb}{1.00,0.94,0.86}
\definecolor{AntiqueWhite2}{rgb}{0.93,0.87,0.80}
\definecolor{AntiqueWhite3}{rgb}{0.80,0.75,0.69}
\definecolor{AntiqueWhite4}{rgb}{0.55,0.51,0.47}
\definecolor{AntiqueWhite}{rgb}{0.98,0.92,0.84}
\definecolor{BlanchedAlmond}{rgb}{1.00,0.92,0.80}
\definecolor{BlueViolet}{rgb}{0.54,0.17,0.89}
\definecolor{CadetBlue1}{rgb}{0.60,0.96,1.00}
\definecolor{CadetBlue2}{rgb}{0.56,0.90,0.93}
\definecolor{CadetBlue3}{rgb}{0.48,0.77,0.80}
\definecolor{CadetBlue4}{rgb}{0.33,0.53,0.55}
\definecolor{CadetBlue}{rgb}{0.37,0.62,0.63}
\definecolor{CornflowerBlue}{rgb}{0.39,0.58,0.93}
\definecolor{DarkBlue}{rgb}{0.00,0.00,0.55}
\definecolor{DarkCyan}{rgb}{0.00,0.55,0.55}
\definecolor{DarkGoldenrod1}{rgb}{1.00,0.73,0.06}
\definecolor{DarkGoldenrod2}{rgb}{0.93,0.68,0.05}
\definecolor{DarkGoldenrod3}{rgb}{0.80,0.58,0.05}
\definecolor{DarkGoldenrod4}{rgb}{0.55,0.40,0.03}
\definecolor{DarkGoldenrod}{rgb}{0.72,0.53,0.04}
\definecolor{DarkGray}{rgb}{0.66,0.66,0.66}
\definecolor{DarkGreen}{rgb}{0.00,0.39,0.00}
\definecolor{DarkGrey}{rgb}{0.66,0.66,0.66}
\definecolor{DarkKhaki}{rgb}{0.74,0.72,0.42}
\definecolor{DarkMagenta}{rgb}{0.55,0.00,0.55}
\definecolor{DarkOliveGreen1}{rgb}{0.79,1.00,0.44}
\definecolor{DarkOliveGreen2}{rgb}{0.74,0.93,0.41}
\definecolor{DarkOliveGreen3}{rgb}{0.64,0.80,0.35}
\definecolor{DarkOliveGreen4}{rgb}{0.43,0.55,0.24}
\definecolor{DarkOliveGreen}{rgb}{0.33,0.42,0.18}
\definecolor{DarkOrange1}{rgb}{1.00,0.50,0.00}
\definecolor{DarkOrange2}{rgb}{0.93,0.46,0.00}
\definecolor{DarkOrange3}{rgb}{0.80,0.40,0.00}
\definecolor{DarkOrange4}{rgb}{0.55,0.27,0.00}
\definecolor{DarkOrange}{rgb}{1.00,0.55,0.00}
\definecolor{DarkOrchid1}{rgb}{0.75,0.24,1.00}
\definecolor{DarkOrchid2}{rgb}{0.70,0.23,0.93}
\definecolor{DarkOrchid3}{rgb}{0.60,0.20,0.80}
\definecolor{DarkOrchid4}{rgb}{0.41,0.13,0.55}
\definecolor{DarkOrchid}{rgb}{0.60,0.20,0.80}
\definecolor{DarkRed}{rgb}{0.55,0.00,0.00}
\definecolor{DarkSalmon}{rgb}{0.91,0.59,0.48}
\definecolor{DarkSeaGreen1}{rgb}{0.76,1.00,0.76}
\definecolor{DarkSeaGreen2}{rgb}{0.71,0.93,0.71}
\definecolor{DarkSeaGreen3}{rgb}{0.61,0.80,0.61}
\definecolor{DarkSeaGreen4}{rgb}{0.41,0.55,0.41}
\definecolor{DarkSeaGreen}{rgb}{0.56,0.74,0.56}
\definecolor{DarkSlateBlue}{rgb}{0.28,0.24,0.55}
\definecolor{DarkSlateGray1}{rgb}{0.59,1.00,1.00}
\definecolor{DarkSlateGray2}{rgb}{0.55,0.93,0.93}
\definecolor{DarkSlateGray3}{rgb}{0.47,0.80,0.80}
\definecolor{DarkSlateGray4}{rgb}{0.32,0.55,0.55}
\definecolor{DarkSlateGray}{rgb}{0.18,0.31,0.31}
\definecolor{DarkSlateGrey}{rgb}{0.18,0.31,0.31}
\definecolor{DarkTurquoise}{rgb}{0.00,0.81,0.82}
\definecolor{DarkViolet}{rgb}{0.58,0.00,0.83}
\definecolor{DeepPink1}{rgb}{1.00,0.08,0.58}
\definecolor{DeepPink2}{rgb}{0.93,0.07,0.54}
\definecolor{DeepPink3}{rgb}{0.80,0.06,0.46}
\definecolor{DeepPink4}{rgb}{0.55,0.04,0.31}
\definecolor{DeepPink}{rgb}{1.00,0.08,0.58}
\definecolor{DeepSkyBlue1}{rgb}{0.00,0.75,1.00}
\definecolor{DeepSkyBlue2}{rgb}{0.00,0.70,0.93}
\definecolor{DeepSkyBlue3}{rgb}{0.00,0.60,0.80}
\definecolor{DeepSkyBlue4}{rgb}{0.00,0.41,0.55}
\definecolor{DeepSkyBlue}{rgb}{0.00,0.75,1.00}
\definecolor{DimGray}{rgb}{0.41,0.41,0.41}
\definecolor{DimGrey}{rgb}{0.41,0.41,0.41}
\definecolor{DodgerBlue1}{rgb}{0.12,0.56,1.00}
\definecolor{DodgerBlue2}{rgb}{0.11,0.53,0.93}
\definecolor{DodgerBlue3}{rgb}{0.09,0.45,0.80}
\definecolor{DodgerBlue4}{rgb}{0.06,0.31,0.55}
\definecolor{DodgerBlue}{rgb}{0.12,0.56,1.00}
\definecolor{FloralWhite}{rgb}{1.00,0.98,0.94}
\definecolor{ForestGreen}{rgb}{0.13,0.55,0.13}
\definecolor{GhostWhite}{rgb}{0.97,0.97,1.00}
\definecolor{GreenYellow}{rgb}{0.68,1.00,0.18}
\definecolor{HotPink1}{rgb}{1.00,0.43,0.71}
\definecolor{HotPink2}{rgb}{0.93,0.42,0.65}
\definecolor{HotPink3}{rgb}{0.80,0.38,0.56}
\definecolor{HotPink4}{rgb}{0.55,0.23,0.38}
\definecolor{HotPink}{rgb}{1.00,0.41,0.71}
\definecolor{IndianRed1}{rgb}{1.00,0.42,0.42}
\definecolor{IndianRed2}{rgb}{0.93,0.39,0.39}
\definecolor{IndianRed3}{rgb}{0.80,0.33,0.33}
\definecolor{IndianRed4}{rgb}{0.55,0.23,0.23}
\definecolor{IndianRed}{rgb}{0.80,0.36,0.36}
\definecolor{LavenderBlush1}{rgb}{1.00,0.94,0.96}
\definecolor{LavenderBlush2}{rgb}{0.93,0.88,0.90}
\definecolor{LavenderBlush3}{rgb}{0.80,0.76,0.77}
\definecolor{LavenderBlush4}{rgb}{0.55,0.51,0.53}
\definecolor{LavenderBlush}{rgb}{1.00,0.94,0.96}
\definecolor{LawnGreen}{rgb}{0.49,0.99,0.00}
\definecolor{LemonChiffon1}{rgb}{1.00,0.98,0.80}
\definecolor{LemonChiffon2}{rgb}{0.93,0.91,0.75}
\definecolor{LemonChiffon3}{rgb}{0.80,0.79,0.65}
\definecolor{LemonChiffon4}{rgb}{0.55,0.54,0.44}
\definecolor{LemonChiffon}{rgb}{1.00,0.98,0.80}
\definecolor{LightBlue1}{rgb}{0.75,0.94,1.00}
\definecolor{LightBlue2}{rgb}{0.70,0.87,0.93}
\definecolor{LightBlue3}{rgb}{0.60,0.75,0.80}
\definecolor{LightBlue4}{rgb}{0.41,0.51,0.55}
\definecolor{LightBlue}{rgb}{0.68,0.85,0.90}
\definecolor{LightCoral}{rgb}{0.94,0.50,0.50}
\definecolor{LightCyan1}{rgb}{0.88,1.00,1.00}
\definecolor{LightCyan2}{rgb}{0.82,0.93,0.93}
\definecolor{LightCyan3}{rgb}{0.71,0.80,0.80}
\definecolor{LightCyan4}{rgb}{0.48,0.55,0.55}
\definecolor{LightCyan}{rgb}{0.88,1.00,1.00}
\definecolor{LightGoldenrod1}{rgb}{1.00,0.93,0.55}
\definecolor{LightGoldenrod2}{rgb}{0.93,0.86,0.51}
\definecolor{LightGoldenrod3}{rgb}{0.80,0.75,0.44}
\definecolor{LightGoldenrod4}{rgb}{0.55,0.51,0.30}
\definecolor{LightGoldenrodYellow}{rgb}{0.98,0.98,0.82}
\definecolor{LightGoldenrod}{rgb}{0.93,0.87,0.51}
\definecolor{LightGray}{rgb}{0.83,0.83,0.83}
\definecolor{LightGreen}{rgb}{0.56,0.93,0.56}
\definecolor{LightGrey}{rgb}{0.83,0.83,0.83}
\definecolor{LightPink1}{rgb}{1.00,0.68,0.73}
\definecolor{LightPink2}{rgb}{0.93,0.64,0.68}
\definecolor{LightPink3}{rgb}{0.80,0.55,0.58}
\definecolor{LightPink4}{rgb}{0.55,0.37,0.40}
\definecolor{LightPink}{rgb}{1.00,0.71,0.76}
\definecolor{LightSalmon1}{rgb}{1.00,0.63,0.48}
\definecolor{LightSalmon2}{rgb}{0.93,0.58,0.45}
\definecolor{LightSalmon3}{rgb}{0.80,0.51,0.38}
\definecolor{LightSalmon4}{rgb}{0.55,0.34,0.26}
\definecolor{LightSalmon}{rgb}{1.00,0.63,0.48}
\definecolor{LightSeaGreen}{rgb}{0.13,0.70,0.67}
\definecolor{LightSkyBlue1}{rgb}{0.69,0.89,1.00}
\definecolor{LightSkyBlue2}{rgb}{0.64,0.83,0.93}
\definecolor{LightSkyBlue3}{rgb}{0.55,0.71,0.80}
\definecolor{LightSkyBlue4}{rgb}{0.38,0.48,0.55}
\definecolor{LightSkyBlue}{rgb}{0.53,0.81,0.98}
\definecolor{LightSlateBlue}{rgb}{0.52,0.44,1.00}
\definecolor{LightSlateGray}{rgb}{0.47,0.53,0.60}
\definecolor{LightSlateGrey}{rgb}{0.47,0.53,0.60}
\definecolor{LightSteelBlue1}{rgb}{0.79,0.88,1.00}
\definecolor{LightSteelBlue2}{rgb}{0.74,0.82,0.93}
\definecolor{LightSteelBlue3}{rgb}{0.64,0.71,0.80}
\definecolor{LightSteelBlue4}{rgb}{0.43,0.48,0.55}
\definecolor{LightSteelBlue}{rgb}{0.69,0.77,0.87}
\definecolor{LightYellow1}{rgb}{1.00,1.00,0.88}
\definecolor{LightYellow2}{rgb}{0.93,0.93,0.82}
\definecolor{LightYellow3}{rgb}{0.80,0.80,0.71}
\definecolor{LightYellow4}{rgb}{0.55,0.55,0.48}
\definecolor{LightYellow}{rgb}{1.00,1.00,0.88}
\definecolor{LimeGreen}{rgb}{0.20,0.80,0.20}
\definecolor{MediumAquamarine}{rgb}{0.40,0.80,0.67}
\definecolor{MediumBlue}{rgb}{0.00,0.00,0.80}
\definecolor{MediumOrchid1}{rgb}{0.88,0.40,1.00}
\definecolor{MediumOrchid2}{rgb}{0.82,0.37,0.93}
\definecolor{MediumOrchid3}{rgb}{0.71,0.32,0.80}
\definecolor{MediumOrchid4}{rgb}{0.48,0.22,0.55}
\definecolor{MediumOrchid}{rgb}{0.73,0.33,0.83}
\definecolor{MediumPurple1}{rgb}{0.67,0.51,1.00}
\definecolor{MediumPurple2}{rgb}{0.62,0.47,0.93}
\definecolor{MediumPurple3}{rgb}{0.54,0.41,0.80}
\definecolor{MediumPurple4}{rgb}{0.36,0.28,0.55}
\definecolor{MediumPurple}{rgb}{0.58,0.44,0.86}
\definecolor{MediumSeaGreen}{rgb}{0.24,0.70,0.44}
\definecolor{MediumSlateBlue}{rgb}{0.48,0.41,0.93}
\definecolor{MediumSpringGreen}{rgb}{0.00,0.98,0.60}
\definecolor{MediumTurquoise}{rgb}{0.28,0.82,0.80}
\definecolor{MediumVioletRed}{rgb}{0.78,0.08,0.52}
\definecolor{MidnightBlue}{rgb}{0.10,0.10,0.44}
\definecolor{MintCream}{rgb}{0.96,1.00,0.98}
\definecolor{MistyRose1}{rgb}{1.00,0.89,0.88}
\definecolor{MistyRose2}{rgb}{0.93,0.84,0.82}
\definecolor{MistyRose3}{rgb}{0.80,0.72,0.71}
\definecolor{MistyRose4}{rgb}{0.55,0.49,0.48}
\definecolor{MistyRose}{rgb}{1.00,0.89,0.88}
\definecolor{NavajoWhite1}{rgb}{1.00,0.87,0.68}
\definecolor{NavajoWhite2}{rgb}{0.93,0.81,0.63}
\definecolor{NavajoWhite3}{rgb}{0.80,0.70,0.55}
\definecolor{NavajoWhite4}{rgb}{0.55,0.47,0.37}
\definecolor{NavajoWhite}{rgb}{1.00,0.87,0.68}
\definecolor{NavyBlue}{rgb}{0.00,0.00,0.50}
\definecolor{OldLace}{rgb}{0.99,0.96,0.90}
\definecolor{OliveDrab1}{rgb}{0.75,1.00,0.24}
\definecolor{OliveDrab2}{rgb}{0.70,0.93,0.23}
\definecolor{OliveDrab3}{rgb}{0.60,0.80,0.20}
\definecolor{OliveDrab4}{rgb}{0.41,0.55,0.13}
\definecolor{OliveDrab}{rgb}{0.42,0.56,0.14}
\definecolor{OrangeRed1}{rgb}{1.00,0.27,0.00}
\definecolor{OrangeRed2}{rgb}{0.93,0.25,0.00}
\definecolor{OrangeRed3}{rgb}{0.80,0.22,0.00}
\definecolor{OrangeRed4}{rgb}{0.55,0.15,0.00}
\definecolor{OrangeRed}{rgb}{1.00,0.27,0.00}
\definecolor{PaleGoldenrod}{rgb}{0.93,0.91,0.67}
\definecolor{PaleGreen1}{rgb}{0.60,1.00,0.60}
\definecolor{PaleGreen2}{rgb}{0.56,0.93,0.56}
\definecolor{PaleGreen3}{rgb}{0.49,0.80,0.49}
\definecolor{PaleGreen4}{rgb}{0.33,0.55,0.33}
\definecolor{PaleGreen}{rgb}{0.60,0.98,0.60}
\definecolor{PaleTurquoise1}{rgb}{0.73,1.00,1.00}
\definecolor{PaleTurquoise2}{rgb}{0.68,0.93,0.93}
\definecolor{PaleTurquoise3}{rgb}{0.59,0.80,0.80}
\definecolor{PaleTurquoise4}{rgb}{0.40,0.55,0.55}
\definecolor{PaleTurquoise}{rgb}{0.69,0.93,0.93}
\definecolor{PaleVioletRed1}{rgb}{1.00,0.51,0.67}
\definecolor{PaleVioletRed2}{rgb}{0.93,0.47,0.62}
\definecolor{PaleVioletRed3}{rgb}{0.80,0.41,0.54}
\definecolor{PaleVioletRed4}{rgb}{0.55,0.28,0.36}
\definecolor{PaleVioletRed}{rgb}{0.86,0.44,0.58}
\definecolor{PapayaWhip}{rgb}{1.00,0.94,0.84}
\definecolor{PeachPuff1}{rgb}{1.00,0.85,0.73}
\definecolor{PeachPuff2}{rgb}{0.93,0.80,0.68}
\definecolor{PeachPuff3}{rgb}{0.80,0.69,0.58}
\definecolor{PeachPuff4}{rgb}{0.55,0.47,0.40}
\definecolor{PeachPuff}{rgb}{1.00,0.85,0.73}
\definecolor{PowderBlue}{rgb}{0.69,0.88,0.90}
\definecolor{RosyBrown1}{rgb}{1.00,0.76,0.76}
\definecolor{RosyBrown2}{rgb}{0.93,0.71,0.71}
\definecolor{RosyBrown3}{rgb}{0.80,0.61,0.61}
\definecolor{RosyBrown4}{rgb}{0.55,0.41,0.41}
\definecolor{RosyBrown}{rgb}{0.74,0.56,0.56}
\definecolor{RoyalBlue1}{rgb}{0.28,0.46,1.00}
\definecolor{RoyalBlue2}{rgb}{0.26,0.43,0.93}
\definecolor{RoyalBlue3}{rgb}{0.23,0.37,0.80}
\definecolor{RoyalBlue4}{rgb}{0.15,0.25,0.55}
\definecolor{RoyalBlue}{rgb}{0.25,0.41,0.88}
\definecolor{SaddleBrown}{rgb}{0.55,0.27,0.07}
\definecolor{SandyBrown}{rgb}{0.96,0.64,0.38}
\definecolor{SeaGreen1}{rgb}{0.33,1.00,0.62}
\definecolor{SeaGreen2}{rgb}{0.31,0.93,0.58}
\definecolor{SeaGreen3}{rgb}{0.26,0.80,0.50}
\definecolor{SeaGreen4}{rgb}{0.18,0.55,0.34}
\definecolor{SeaGreen}{rgb}{0.18,0.55,0.34}
\definecolor{SkyBlue1}{rgb}{0.53,0.81,1.00}
\definecolor{SkyBlue2}{rgb}{0.49,0.75,0.93}
\definecolor{SkyBlue3}{rgb}{0.42,0.65,0.80}
\definecolor{SkyBlue4}{rgb}{0.29,0.44,0.55}
\definecolor{SkyBlue}{rgb}{0.53,0.81,0.92}
\definecolor{SlateBlue1}{rgb}{0.51,0.44,1.00}
\definecolor{SlateBlue2}{rgb}{0.48,0.40,0.93}
\definecolor{SlateBlue3}{rgb}{0.41,0.35,0.80}
\definecolor{SlateBlue4}{rgb}{0.28,0.24,0.55}
\definecolor{SlateBlue}{rgb}{0.42,0.35,0.80}
\definecolor{SlateGray1}{rgb}{0.78,0.89,1.00}
\definecolor{SlateGray2}{rgb}{0.73,0.83,0.93}
\definecolor{SlateGray3}{rgb}{0.62,0.71,0.80}
\definecolor{SlateGray4}{rgb}{0.42,0.48,0.55}
\definecolor{SlateGray}{rgb}{0.44,0.50,0.56}
\definecolor{SlateGrey}{rgb}{0.44,0.50,0.56}
\definecolor{SpringGreen1}{rgb}{0.00,1.00,0.50}
\definecolor{SpringGreen2}{rgb}{0.00,0.93,0.46}
\definecolor{SpringGreen3}{rgb}{0.00,0.80,0.40}
\definecolor{SpringGreen4}{rgb}{0.00,0.55,0.27}
\definecolor{SpringGreen}{rgb}{0.00,1.00,0.50}
\definecolor{SteelBlue1}{rgb}{0.39,0.72,1.00}
\definecolor{SteelBlue2}{rgb}{0.36,0.67,0.93}
\definecolor{SteelBlue3}{rgb}{0.31,0.58,0.80}
\definecolor{SteelBlue4}{rgb}{0.21,0.39,0.55}
\definecolor{SteelBlue}{rgb}{0.27,0.51,0.71}
\definecolor{VioletRed1}{rgb}{1.00,0.24,0.59}
\definecolor{VioletRed2}{rgb}{0.93,0.23,0.55}
\definecolor{VioletRed3}{rgb}{0.80,0.20,0.47}
\definecolor{VioletRed4}{rgb}{0.55,0.13,0.32}
\definecolor{VioletRed}{rgb}{0.82,0.13,0.56}
\definecolor{WhiteSmoke}{rgb}{0.96,0.96,0.96}
\definecolor{YellowGreen}{rgb}{0.60,0.80,0.20}
\definecolor{aliceblue}{rgb}{0.94,0.97,1.00}
\definecolor{antiquewhite}{rgb}{0.98,0.92,0.84}
\definecolor{aquamarine1}{rgb}{0.50,1.00,0.83}
\definecolor{aquamarine2}{rgb}{0.46,0.93,0.78}
\definecolor{aquamarine3}{rgb}{0.40,0.80,0.67}
\definecolor{aquamarine4}{rgb}{0.27,0.55,0.45}
\definecolor{aquamarine}{rgb}{0.50,1.00,0.83}
\definecolor{azure1}{rgb}{0.94,1.00,1.00}
\definecolor{azure2}{rgb}{0.88,0.93,0.93}
\definecolor{azure3}{rgb}{0.76,0.80,0.80}
\definecolor{azure4}{rgb}{0.51,0.55,0.55}
\definecolor{azure}{rgb}{0.94,1.00,1.00}
\definecolor{beige}{rgb}{0.96,0.96,0.86}
\definecolor{bisque1}{rgb}{1.00,0.89,0.77}
\definecolor{bisque2}{rgb}{0.93,0.84,0.72}
\definecolor{bisque3}{rgb}{0.80,0.72,0.62}
\definecolor{bisque4}{rgb}{0.55,0.49,0.42}
\definecolor{bisque}{rgb}{1.00,0.89,0.77}
\definecolor{black}{rgb}{0.00,0.00,0.00}
\definecolor{blanchedalmond}{rgb}{1.00,0.92,0.80}
\definecolor{blue1}{rgb}{0.00,0.00,1.00}
\definecolor{blue2}{rgb}{0.00,0.00,0.93}
\definecolor{blue3}{rgb}{0.00,0.00,0.80}
\definecolor{blue4}{rgb}{0.00,0.00,0.55}
\definecolor{blueviolet}{rgb}{0.54,0.17,0.89}
\definecolor{blue}{rgb}{0.00,0.00,1.00}
\definecolor{brown1}{rgb}{1.00,0.25,0.25}
\definecolor{brown2}{rgb}{0.93,0.23,0.23}
\definecolor{brown3}{rgb}{0.80,0.20,0.20}
\definecolor{brown4}{rgb}{0.55,0.14,0.14}
\definecolor{brown}{rgb}{0.65,0.16,0.16}
\definecolor{burlywood1}{rgb}{1.00,0.83,0.61}
\definecolor{burlywood2}{rgb}{0.93,0.77,0.57}
\definecolor{burlywood3}{rgb}{0.80,0.67,0.49}
\definecolor{burlywood4}{rgb}{0.55,0.45,0.33}
\definecolor{burlywood}{rgb}{0.87,0.72,0.53}
\definecolor{cadetblue}{rgb}{0.37,0.62,0.63}
\definecolor{chartreuse1}{rgb}{0.50,1.00,0.00}
\definecolor{chartreuse2}{rgb}{0.46,0.93,0.00}
\definecolor{chartreuse3}{rgb}{0.40,0.80,0.00}
\definecolor{chartreuse4}{rgb}{0.27,0.55,0.00}
\definecolor{chartreuse}{rgb}{0.50,1.00,0.00}
\definecolor{chocolate1}{rgb}{1.00,0.50,0.14}
\definecolor{chocolate2}{rgb}{0.93,0.46,0.13}
\definecolor{chocolate3}{rgb}{0.80,0.40,0.11}
\definecolor{chocolate4}{rgb}{0.55,0.27,0.07}
\definecolor{chocolate}{rgb}{0.82,0.41,0.12}
\definecolor{coral1}{rgb}{1.00,0.45,0.34}
\definecolor{coral2}{rgb}{0.93,0.42,0.31}
\definecolor{coral3}{rgb}{0.80,0.36,0.27}
\definecolor{coral4}{rgb}{0.55,0.24,0.18}
\definecolor{coral}{rgb}{1.00,0.50,0.31}
\definecolor{cornflowerblue}{rgb}{0.39,0.58,0.93}
\definecolor{cornsilk1}{rgb}{1.00,0.97,0.86}
\definecolor{cornsilk2}{rgb}{0.93,0.91,0.80}
\definecolor{cornsilk3}{rgb}{0.80,0.78,0.69}
\definecolor{cornsilk4}{rgb}{0.55,0.53,0.47}
\definecolor{cornsilk}{rgb}{1.00,0.97,0.86}
\definecolor{cyan1}{rgb}{0.00,1.00,1.00}
\definecolor{cyan2}{rgb}{0.00,0.93,0.93}
\definecolor{cyan3}{rgb}{0.00,0.80,0.80}
\definecolor{cyan4}{rgb}{0.00,0.55,0.55}
\definecolor{cyan}{rgb}{0.00,1.00,1.00}
\definecolor{darkblue}{rgb}{0.00,0.00,0.55}
\definecolor{darkcyan}{rgb}{0.00,0.55,0.55}
\definecolor{darkgoldenrod}{rgb}{0.72,0.53,0.04}
\definecolor{darkgray}{rgb}{0.66,0.66,0.66}
\definecolor{darkgreen}{rgb}{0.00,0.39,0.00}
\definecolor{darkgrey}{rgb}{0.66,0.66,0.66}
\definecolor{darkkhaki}{rgb}{0.74,0.72,0.42}
\definecolor{darkmagenta}{rgb}{0.55,0.00,0.55}
\definecolor{darkolive}{rgb}{0.33,0.42,0.18}
\definecolor{darkorange}{rgb}{1.00,0.55,0.00}
\definecolor{darkorchid}{rgb}{0.60,0.20,0.80}
\definecolor{darkred}{rgb}{0.55,0.00,0.00}
\definecolor{darksalmon}{rgb}{0.91,0.59,0.48}
\definecolor{darksea}{rgb}{0.56,0.74,0.56}
\definecolor{darkslate}{rgb}{0.18,0.31,0.31}
\definecolor{darkslate}{rgb}{0.18,0.31,0.31}
\definecolor{darkslate}{rgb}{0.28,0.24,0.55}
\definecolor{darkturquoise}{rgb}{0.00,0.81,0.82}
\definecolor{darkviolet}{rgb}{0.58,0.00,0.83}
\definecolor{deeppink}{rgb}{1.00,0.08,0.58}
\definecolor{deepsky}{rgb}{0.00,0.75,1.00}
\definecolor{dimgray}{rgb}{0.41,0.41,0.41}
\definecolor{dimgrey}{rgb}{0.41,0.41,0.41}
\definecolor{dodgerblue}{rgb}{0.12,0.56,1.00}
\definecolor{firebrick1}{rgb}{1.00,0.19,0.19}
\definecolor{firebrick2}{rgb}{0.93,0.17,0.17}
\definecolor{firebrick3}{rgb}{0.80,0.15,0.15}
\definecolor{firebrick4}{rgb}{0.55,0.10,0.10}
\definecolor{firebrick}{rgb}{0.70,0.13,0.13}
\definecolor{floralwhite}{rgb}{1.00,0.98,0.94}
\definecolor{forestgreen}{rgb}{0.13,0.55,0.13}
\definecolor{gainsboro}{rgb}{0.86,0.86,0.86}
\definecolor{ghostwhite}{rgb}{0.97,0.97,1.00}
\definecolor{gold1}{rgb}{1.00,0.84,0.00}
\definecolor{gold2}{rgb}{0.93,0.79,0.00}
\definecolor{gold3}{rgb}{0.80,0.68,0.00}
\definecolor{gold4}{rgb}{0.55,0.46,0.00}
\definecolor{goldenrod1}{rgb}{1.00,0.76,0.15}
\definecolor{goldenrod2}{rgb}{0.93,0.71,0.13}
\definecolor{goldenrod3}{rgb}{0.80,0.61,0.11}
\definecolor{goldenrod4}{rgb}{0.55,0.41,0.08}
\definecolor{goldenrod}{rgb}{0.85,0.65,0.13}
\definecolor{gold}{rgb}{1.00,0.84,0.00}
\definecolor{gray0}{rgb}{0.00,0.00,0.00}
\definecolor{gray100}{rgb}{1.00,1.00,1.00}
\definecolor{gray10}{rgb}{0.10,0.10,0.10}
\definecolor{gray11}{rgb}{0.11,0.11,0.11}
\definecolor{gray12}{rgb}{0.12,0.12,0.12}
\definecolor{gray13}{rgb}{0.13,0.13,0.13}
\definecolor{gray14}{rgb}{0.14,0.14,0.14}
\definecolor{gray15}{rgb}{0.15,0.15,0.15}
\definecolor{gray16}{rgb}{0.16,0.16,0.16}
\definecolor{gray17}{rgb}{0.17,0.17,0.17}
\definecolor{gray18}{rgb}{0.18,0.18,0.18}
\definecolor{gray19}{rgb}{0.19,0.19,0.19}
\definecolor{gray1}{rgb}{0.01,0.01,0.01}
\definecolor{gray20}{rgb}{0.20,0.20,0.20}
\definecolor{gray21}{rgb}{0.21,0.21,0.21}
\definecolor{gray22}{rgb}{0.22,0.22,0.22}
\definecolor{gray23}{rgb}{0.23,0.23,0.23}
\definecolor{gray24}{rgb}{0.24,0.24,0.24}
\definecolor{gray25}{rgb}{0.25,0.25,0.25}
\definecolor{gray26}{rgb}{0.26,0.26,0.26}
\definecolor{gray27}{rgb}{0.27,0.27,0.27}
\definecolor{gray28}{rgb}{0.28,0.28,0.28}
\definecolor{gray29}{rgb}{0.29,0.29,0.29}
\definecolor{gray2}{rgb}{0.02,0.02,0.02}
\definecolor{gray30}{rgb}{0.30,0.30,0.30}
\definecolor{gray31}{rgb}{0.31,0.31,0.31}
\definecolor{gray32}{rgb}{0.32,0.32,0.32}
\definecolor{gray33}{rgb}{0.33,0.33,0.33}
\definecolor{gray34}{rgb}{0.34,0.34,0.34}
\definecolor{gray35}{rgb}{0.35,0.35,0.35}
\definecolor{gray36}{rgb}{0.36,0.36,0.36}
\definecolor{gray37}{rgb}{0.37,0.37,0.37}
\definecolor{gray38}{rgb}{0.38,0.38,0.38}
\definecolor{gray39}{rgb}{0.39,0.39,0.39}
\definecolor{gray3}{rgb}{0.03,0.03,0.03}
\definecolor{gray40}{rgb}{0.40,0.40,0.40}
\definecolor{gray41}{rgb}{0.41,0.41,0.41}
\definecolor{gray42}{rgb}{0.42,0.42,0.42}
\definecolor{gray43}{rgb}{0.43,0.43,0.43}
\definecolor{gray44}{rgb}{0.44,0.44,0.44}
\definecolor{gray45}{rgb}{0.45,0.45,0.45}
\definecolor{gray46}{rgb}{0.46,0.46,0.46}
\definecolor{gray47}{rgb}{0.47,0.47,0.47}
\definecolor{gray48}{rgb}{0.48,0.48,0.48}
\definecolor{gray49}{rgb}{0.49,0.49,0.49}
\definecolor{gray4}{rgb}{0.04,0.04,0.04}
\definecolor{gray50}{rgb}{0.50,0.50,0.50}
\definecolor{gray51}{rgb}{0.51,0.51,0.51}
\definecolor{gray52}{rgb}{0.52,0.52,0.52}
\definecolor{gray53}{rgb}{0.53,0.53,0.53}
\definecolor{gray54}{rgb}{0.54,0.54,0.54}
\definecolor{gray55}{rgb}{0.55,0.55,0.55}
\definecolor{gray56}{rgb}{0.56,0.56,0.56}
\definecolor{gray57}{rgb}{0.57,0.57,0.57}
\definecolor{gray58}{rgb}{0.58,0.58,0.58}
\definecolor{gray59}{rgb}{0.59,0.59,0.59}
\definecolor{gray5}{rgb}{0.05,0.05,0.05}
\definecolor{gray60}{rgb}{0.60,0.60,0.60}
\definecolor{gray61}{rgb}{0.61,0.61,0.61}
\definecolor{gray62}{rgb}{0.62,0.62,0.62}
\definecolor{gray63}{rgb}{0.63,0.63,0.63}
\definecolor{gray64}{rgb}{0.64,0.64,0.64}
\definecolor{gray65}{rgb}{0.65,0.65,0.65}
\definecolor{gray66}{rgb}{0.66,0.66,0.66}
\definecolor{gray67}{rgb}{0.67,0.67,0.67}
\definecolor{gray68}{rgb}{0.68,0.68,0.68}
\definecolor{gray69}{rgb}{0.69,0.69,0.69}
\definecolor{gray6}{rgb}{0.06,0.06,0.06}
\definecolor{gray70}{rgb}{0.70,0.70,0.70}
\definecolor{gray71}{rgb}{0.71,0.71,0.71}
\definecolor{gray72}{rgb}{0.72,0.72,0.72}
\definecolor{gray73}{rgb}{0.73,0.73,0.73}
\definecolor{gray74}{rgb}{0.74,0.74,0.74}
\definecolor{gray75}{rgb}{0.75,0.75,0.75}
\definecolor{gray76}{rgb}{0.76,0.76,0.76}
\definecolor{gray77}{rgb}{0.77,0.77,0.77}
\definecolor{gray78}{rgb}{0.78,0.78,0.78}
\definecolor{gray79}{rgb}{0.79,0.79,0.79}
\definecolor{gray7}{rgb}{0.07,0.07,0.07}
\definecolor{gray80}{rgb}{0.80,0.80,0.80}
\definecolor{gray81}{rgb}{0.81,0.81,0.81}
\definecolor{gray82}{rgb}{0.82,0.82,0.82}
\definecolor{gray83}{rgb}{0.83,0.83,0.83}
\definecolor{gray84}{rgb}{0.84,0.84,0.84}
\definecolor{gray85}{rgb}{0.85,0.85,0.85}
\definecolor{gray86}{rgb}{0.86,0.86,0.86}
\definecolor{gray87}{rgb}{0.87,0.87,0.87}
\definecolor{gray88}{rgb}{0.88,0.88,0.88}
\definecolor{gray89}{rgb}{0.89,0.89,0.89}
\definecolor{gray8}{rgb}{0.08,0.08,0.08}
\definecolor{gray90}{rgb}{0.90,0.90,0.90}
\definecolor{gray91}{rgb}{0.91,0.91,0.91}
\definecolor{gray92}{rgb}{0.92,0.92,0.92}
\definecolor{gray93}{rgb}{0.93,0.93,0.93}
\definecolor{gray94}{rgb}{0.94,0.94,0.94}
\definecolor{gray95}{rgb}{0.95,0.95,0.95}
\definecolor{gray96}{rgb}{0.96,0.96,0.96}
\definecolor{gray97}{rgb}{0.97,0.97,0.97}
\definecolor{gray98}{rgb}{0.98,0.98,0.98}
\definecolor{gray99}{rgb}{0.99,0.99,0.99}
\definecolor{gray9}{rgb}{0.09,0.09,0.09}
\definecolor{gray}{rgb}{0.75,0.75,0.75}
\definecolor{green1}{rgb}{0.00,1.00,0.00}
\definecolor{green2}{rgb}{0.00,0.93,0.00}
\definecolor{green3}{rgb}{0.00,0.80,0.00}
\definecolor{green4}{rgb}{0.00,0.55,0.00}
\definecolor{greenyellow}{rgb}{0.68,1.00,0.18}
\definecolor{green}{rgb}{0.00,1.00,0.00}
\definecolor{grey0}{rgb}{0.00,0.00,0.00}
\definecolor{grey100}{rgb}{1.00,1.00,1.00}
\definecolor{grey10}{rgb}{0.10,0.10,0.10}
\definecolor{grey11}{rgb}{0.11,0.11,0.11}
\definecolor{grey12}{rgb}{0.12,0.12,0.12}
\definecolor{grey13}{rgb}{0.13,0.13,0.13}
\definecolor{grey14}{rgb}{0.14,0.14,0.14}
\definecolor{grey15}{rgb}{0.15,0.15,0.15}
\definecolor{grey16}{rgb}{0.16,0.16,0.16}
\definecolor{grey17}{rgb}{0.17,0.17,0.17}
\definecolor{grey18}{rgb}{0.18,0.18,0.18}
\definecolor{grey19}{rgb}{0.19,0.19,0.19}
\definecolor{grey1}{rgb}{0.01,0.01,0.01}
\definecolor{grey20}{rgb}{0.20,0.20,0.20}
\definecolor{grey21}{rgb}{0.21,0.21,0.21}
\definecolor{grey22}{rgb}{0.22,0.22,0.22}
\definecolor{grey23}{rgb}{0.23,0.23,0.23}
\definecolor{grey24}{rgb}{0.24,0.24,0.24}
\definecolor{grey25}{rgb}{0.25,0.25,0.25}
\definecolor{grey26}{rgb}{0.26,0.26,0.26}
\definecolor{grey27}{rgb}{0.27,0.27,0.27}
\definecolor{grey28}{rgb}{0.28,0.28,0.28}
\definecolor{grey29}{rgb}{0.29,0.29,0.29}
\definecolor{grey2}{rgb}{0.02,0.02,0.02}
\definecolor{grey30}{rgb}{0.30,0.30,0.30}
\definecolor{grey31}{rgb}{0.31,0.31,0.31}
\definecolor{grey32}{rgb}{0.32,0.32,0.32}
\definecolor{grey33}{rgb}{0.33,0.33,0.33}
\definecolor{grey34}{rgb}{0.34,0.34,0.34}
\definecolor{grey35}{rgb}{0.35,0.35,0.35}
\definecolor{grey36}{rgb}{0.36,0.36,0.36}
\definecolor{grey37}{rgb}{0.37,0.37,0.37}
\definecolor{grey38}{rgb}{0.38,0.38,0.38}
\definecolor{grey39}{rgb}{0.39,0.39,0.39}
\definecolor{grey3}{rgb}{0.03,0.03,0.03}
\definecolor{grey40}{rgb}{0.40,0.40,0.40}
\definecolor{grey41}{rgb}{0.41,0.41,0.41}
\definecolor{grey42}{rgb}{0.42,0.42,0.42}
\definecolor{grey43}{rgb}{0.43,0.43,0.43}
\definecolor{grey44}{rgb}{0.44,0.44,0.44}
\definecolor{grey45}{rgb}{0.45,0.45,0.45}
\definecolor{grey46}{rgb}{0.46,0.46,0.46}
\definecolor{grey47}{rgb}{0.47,0.47,0.47}
\definecolor{grey48}{rgb}{0.48,0.48,0.48}
\definecolor{grey49}{rgb}{0.49,0.49,0.49}
\definecolor{grey4}{rgb}{0.04,0.04,0.04}
\definecolor{grey50}{rgb}{0.50,0.50,0.50}
\definecolor{grey51}{rgb}{0.51,0.51,0.51}
\definecolor{grey52}{rgb}{0.52,0.52,0.52}
\definecolor{grey53}{rgb}{0.53,0.53,0.53}
\definecolor{grey54}{rgb}{0.54,0.54,0.54}
\definecolor{grey55}{rgb}{0.55,0.55,0.55}
\definecolor{grey56}{rgb}{0.56,0.56,0.56}
\definecolor{grey57}{rgb}{0.57,0.57,0.57}
\definecolor{grey58}{rgb}{0.58,0.58,0.58}
\definecolor{grey59}{rgb}{0.59,0.59,0.59}
\definecolor{grey5}{rgb}{0.05,0.05,0.05}
\definecolor{grey60}{rgb}{0.60,0.60,0.60}
\definecolor{grey61}{rgb}{0.61,0.61,0.61}
\definecolor{grey62}{rgb}{0.62,0.62,0.62}
\definecolor{grey63}{rgb}{0.63,0.63,0.63}
\definecolor{grey64}{rgb}{0.64,0.64,0.64}
\definecolor{grey65}{rgb}{0.65,0.65,0.65}
\definecolor{grey66}{rgb}{0.66,0.66,0.66}
\definecolor{grey67}{rgb}{0.67,0.67,0.67}
\definecolor{grey68}{rgb}{0.68,0.68,0.68}
\definecolor{grey69}{rgb}{0.69,0.69,0.69}
\definecolor{grey6}{rgb}{0.06,0.06,0.06}
\definecolor{grey70}{rgb}{0.70,0.70,0.70}
\definecolor{grey71}{rgb}{0.71,0.71,0.71}
\definecolor{grey72}{rgb}{0.72,0.72,0.72}
\definecolor{grey73}{rgb}{0.73,0.73,0.73}
\definecolor{grey74}{rgb}{0.74,0.74,0.74}
\definecolor{grey75}{rgb}{0.75,0.75,0.75}
\definecolor{grey76}{rgb}{0.76,0.76,0.76}
\definecolor{grey77}{rgb}{0.77,0.77,0.77}
\definecolor{grey78}{rgb}{0.78,0.78,0.78}
\definecolor{grey79}{rgb}{0.79,0.79,0.79}
\definecolor{grey7}{rgb}{0.07,0.07,0.07}
\definecolor{grey80}{rgb}{0.80,0.80,0.80}
\definecolor{grey81}{rgb}{0.81,0.81,0.81}
\definecolor{grey82}{rgb}{0.82,0.82,0.82}
\definecolor{grey83}{rgb}{0.83,0.83,0.83}
\definecolor{grey84}{rgb}{0.84,0.84,0.84}
\definecolor{grey85}{rgb}{0.85,0.85,0.85}
\definecolor{grey86}{rgb}{0.86,0.86,0.86}
\definecolor{grey87}{rgb}{0.87,0.87,0.87}
\definecolor{grey88}{rgb}{0.88,0.88,0.88}
\definecolor{grey89}{rgb}{0.89,0.89,0.89}
\definecolor{grey8}{rgb}{0.08,0.08,0.08}
\definecolor{grey90}{rgb}{0.90,0.90,0.90}
\definecolor{grey91}{rgb}{0.91,0.91,0.91}
\definecolor{grey92}{rgb}{0.92,0.92,0.92}
\definecolor{grey93}{rgb}{0.93,0.93,0.93}
\definecolor{grey94}{rgb}{0.94,0.94,0.94}
\definecolor{grey95}{rgb}{0.95,0.95,0.95}
\definecolor{grey96}{rgb}{0.96,0.96,0.96}
\definecolor{grey97}{rgb}{0.97,0.97,0.97}
\definecolor{grey98}{rgb}{0.98,0.98,0.98}
\definecolor{grey99}{rgb}{0.99,0.99,0.99}
\definecolor{grey9}{rgb}{0.09,0.09,0.09}
\definecolor{grey}{rgb}{0.75,0.75,0.75}
\definecolor{honeydew1}{rgb}{0.94,1.00,0.94}
\definecolor{honeydew2}{rgb}{0.88,0.93,0.88}
\definecolor{honeydew3}{rgb}{0.76,0.80,0.76}
\definecolor{honeydew4}{rgb}{0.51,0.55,0.51}
\definecolor{honeydew}{rgb}{0.94,1.00,0.94}
\definecolor{hotpink}{rgb}{1.00,0.41,0.71}
\definecolor{indianred}{rgb}{0.80,0.36,0.36}
\definecolor{ivory1}{rgb}{1.00,1.00,0.94}
\definecolor{ivory2}{rgb}{0.93,0.93,0.88}
\definecolor{ivory3}{rgb}{0.80,0.80,0.76}
\definecolor{ivory4}{rgb}{0.55,0.55,0.51}
\definecolor{ivory}{rgb}{1.00,1.00,0.94}
\definecolor{khaki1}{rgb}{1.00,0.96,0.56}
\definecolor{khaki2}{rgb}{0.93,0.90,0.52}
\definecolor{khaki3}{rgb}{0.80,0.78,0.45}
\definecolor{khaki4}{rgb}{0.55,0.53,0.31}
\definecolor{khaki}{rgb}{0.94,0.90,0.55}
\definecolor{lavenderblush}{rgb}{1.00,0.94,0.96}
\definecolor{lavender}{rgb}{0.90,0.90,0.98}
\definecolor{lawngreen}{rgb}{0.49,0.99,0.00}
\definecolor{lemonchiffon}{rgb}{1.00,0.98,0.80}
\definecolor{lightblue}{rgb}{0.68,0.85,0.90}
\definecolor{lightcoral}{rgb}{0.94,0.50,0.50}
\definecolor{lightcyan}{rgb}{0.88,1.00,1.00}
\definecolor{lightgoldenrod}{rgb}{0.93,0.87,0.51}
\definecolor{lightgoldenrod}{rgb}{0.98,0.98,0.82}
\definecolor{lightgray}{rgb}{0.83,0.83,0.83}
\definecolor{lightgreen}{rgb}{0.56,0.93,0.56}
\definecolor{lightgrey}{rgb}{0.83,0.83,0.83}
\definecolor{lightpink}{rgb}{1.00,0.71,0.76}
\definecolor{lightsalmon}{rgb}{1.00,0.63,0.48}
\definecolor{lightsea}{rgb}{0.13,0.70,0.67}
\definecolor{lightsky}{rgb}{0.53,0.81,0.98}
\definecolor{lightslate}{rgb}{0.47,0.53,0.60}
\definecolor{lightslate}{rgb}{0.47,0.53,0.60}
\definecolor{lightslate}{rgb}{0.52,0.44,1.00}
\definecolor{lightsteel}{rgb}{0.69,0.77,0.87}
\definecolor{lightyellow}{rgb}{1.00,1.00,0.88}
\definecolor{limegreen}{rgb}{0.20,0.80,0.20}
\definecolor{linen}{rgb}{0.98,0.94,0.90}
\definecolor{magenta1}{rgb}{1.00,0.00,1.00}
\definecolor{magenta2}{rgb}{0.93,0.00,0.93}
\definecolor{magenta3}{rgb}{0.80,0.00,0.80}
\definecolor{magenta4}{rgb}{0.55,0.00,0.55}
\definecolor{magenta}{rgb}{1.00,0.00,1.00}
\definecolor{maroon1}{rgb}{1.00,0.20,0.70}
\definecolor{maroon2}{rgb}{0.93,0.19,0.65}
\definecolor{maroon3}{rgb}{0.80,0.16,0.56}
\definecolor{maroon4}{rgb}{0.55,0.11,0.38}
\definecolor{maroon}{rgb}{0.69,0.19,0.38}
\definecolor{mediumaquamarine}{rgb}{0.40,0.80,0.67}
\definecolor{mediumblue}{rgb}{0.00,0.00,0.80}
\definecolor{mediumorchid}{rgb}{0.73,0.33,0.83}
\definecolor{mediumpurple}{rgb}{0.58,0.44,0.86}
\definecolor{mediumsea}{rgb}{0.24,0.70,0.44}
\definecolor{mediumslate}{rgb}{0.48,0.41,0.93}
\definecolor{mediumspring}{rgb}{0.00,0.98,0.60}
\definecolor{mediumturquoise}{rgb}{0.28,0.82,0.80}
\definecolor{mediumviolet}{rgb}{0.78,0.08,0.52}
\definecolor{midnightblue}{rgb}{0.10,0.10,0.44}
\definecolor{mintcream}{rgb}{0.96,1.00,0.98}
\definecolor{mistyrose}{rgb}{1.00,0.89,0.88}
\definecolor{moccasin}{rgb}{1.00,0.89,0.71}
\definecolor{navajowhite}{rgb}{1.00,0.87,0.68}
\definecolor{navyblue}{rgb}{0.00,0.00,0.50}
\definecolor{navy}{rgb}{0.00,0.00,0.50}
\definecolor{oldlace}{rgb}{0.99,0.96,0.90}
\definecolor{olivedrab}{rgb}{0.42,0.56,0.14}
\definecolor{orange1}{rgb}{1.00,0.65,0.00}
\definecolor{orange2}{rgb}{0.93,0.60,0.00}
\definecolor{orange3}{rgb}{0.80,0.52,0.00}
\definecolor{orange4}{rgb}{0.55,0.35,0.00}
\definecolor{orangered}{rgb}{1.00,0.27,0.00}
\definecolor{orange}{rgb}{1.00,0.65,0.00}
\definecolor{orchid1}{rgb}{1.00,0.51,0.98}
\definecolor{orchid2}{rgb}{0.93,0.48,0.91}
\definecolor{orchid3}{rgb}{0.80,0.41,0.79}
\definecolor{orchid4}{rgb}{0.55,0.28,0.54}
\definecolor{orchid}{rgb}{0.85,0.44,0.84}
\definecolor{palegoldenrod}{rgb}{0.93,0.91,0.67}
\definecolor{palegreen}{rgb}{0.60,0.98,0.60}
\definecolor{paleturquoise}{rgb}{0.69,0.93,0.93}
\definecolor{paleviolet}{rgb}{0.86,0.44,0.58}
\definecolor{papayawhip}{rgb}{1.00,0.94,0.84}
\definecolor{peachpuff}{rgb}{1.00,0.85,0.73}
\definecolor{peru}{rgb}{0.80,0.52,0.25}
\definecolor{pink1}{rgb}{1.00,0.71,0.77}
\definecolor{pink2}{rgb}{0.93,0.66,0.72}
\definecolor{pink3}{rgb}{0.80,0.57,0.62}
\definecolor{pink4}{rgb}{0.55,0.39,0.42}
\definecolor{pink}{rgb}{1.00,0.75,0.80}
\definecolor{plum1}{rgb}{1.00,0.73,1.00}
\definecolor{plum2}{rgb}{0.93,0.68,0.93}
\definecolor{plum3}{rgb}{0.80,0.59,0.80}
\definecolor{plum4}{rgb}{0.55,0.40,0.55}
\definecolor{plum}{rgb}{0.87,0.63,0.87}
\definecolor{powderblue}{rgb}{0.69,0.88,0.90}
\definecolor{purple1}{rgb}{0.61,0.19,1.00}
\definecolor{purple2}{rgb}{0.57,0.17,0.93}
\definecolor{purple3}{rgb}{0.49,0.15,0.80}
\definecolor{purple4}{rgb}{0.33,0.10,0.55}
\definecolor{purple}{rgb}{0.63,0.13,0.94}
\definecolor{red1}{rgb}{1.00,0.00,0.00}
\definecolor{red2}{rgb}{0.93,0.00,0.00}
\definecolor{red3}{rgb}{0.80,0.00,0.00}
\definecolor{red4}{rgb}{0.55,0.00,0.00}
\definecolor{red}{rgb}{1.00,0.00,0.00}
\definecolor{rosybrown}{rgb}{0.74,0.56,0.56}
\definecolor{royalblue}{rgb}{0.25,0.41,0.88}
\definecolor{saddlebrown}{rgb}{0.55,0.27,0.07}
\definecolor{salmon1}{rgb}{1.00,0.55,0.41}
\definecolor{salmon2}{rgb}{0.93,0.51,0.38}
\definecolor{salmon3}{rgb}{0.80,0.44,0.33}
\definecolor{salmon4}{rgb}{0.55,0.30,0.22}
\definecolor{salmon}{rgb}{0.98,0.50,0.45}
\definecolor{sandybrown}{rgb}{0.96,0.64,0.38}
\definecolor{seagreen}{rgb}{0.18,0.55,0.34}
\definecolor{seashell1}{rgb}{1.00,0.96,0.93}
\definecolor{seashell2}{rgb}{0.93,0.90,0.87}
\definecolor{seashell3}{rgb}{0.80,0.77,0.75}
\definecolor{seashell4}{rgb}{0.55,0.53,0.51}
\definecolor{seashell}{rgb}{1.00,0.96,0.93}
\definecolor{sienna1}{rgb}{1.00,0.51,0.28}
\definecolor{sienna2}{rgb}{0.93,0.47,0.26}
\definecolor{sienna3}{rgb}{0.80,0.41,0.22}
\definecolor{sienna4}{rgb}{0.55,0.28,0.15}
\definecolor{sienna}{rgb}{0.63,0.32,0.18}
\definecolor{skyblue}{rgb}{0.53,0.81,0.92}
\definecolor{slateblue}{rgb}{0.42,0.35,0.80}
\definecolor{slategray}{rgb}{0.44,0.50,0.56}
\definecolor{slategrey}{rgb}{0.44,0.50,0.56}
\definecolor{snow1}{rgb}{1.00,0.98,0.98}
\definecolor{snow2}{rgb}{0.93,0.91,0.91}
\definecolor{snow3}{rgb}{0.80,0.79,0.79}
\definecolor{snow4}{rgb}{0.55,0.54,0.54}
\definecolor{snow}{rgb}{1.00,0.98,0.98}
\definecolor{springgreen}{rgb}{0.00,1.00,0.50}
\definecolor{steelblue}{rgb}{0.27,0.51,0.71}
\definecolor{tan1}{rgb}{1.00,0.65,0.31}
\definecolor{tan2}{rgb}{0.93,0.60,0.29}
\definecolor{tan3}{rgb}{0.80,0.52,0.25}
\definecolor{tan4}{rgb}{0.55,0.35,0.17}
\definecolor{tan}{rgb}{0.82,0.71,0.55}
\definecolor{thistle1}{rgb}{1.00,0.88,1.00}
\definecolor{thistle2}{rgb}{0.93,0.82,0.93}
\definecolor{thistle3}{rgb}{0.80,0.71,0.80}
\definecolor{thistle4}{rgb}{0.55,0.48,0.55}
\definecolor{thistle}{rgb}{0.85,0.75,0.85}
\definecolor{tomato1}{rgb}{1.00,0.39,0.28}
\definecolor{tomato2}{rgb}{0.93,0.36,0.26}
\definecolor{tomato3}{rgb}{0.80,0.31,0.22}
\definecolor{tomato4}{rgb}{0.55,0.21,0.15}
\definecolor{tomato}{rgb}{1.00,0.39,0.28}
\definecolor{turquoise1}{rgb}{0.00,0.96,1.00}
\definecolor{turquoise2}{rgb}{0.00,0.90,0.93}
\definecolor{turquoise3}{rgb}{0.00,0.77,0.80}
\definecolor{turquoise4}{rgb}{0.00,0.53,0.55}
\definecolor{turquoise}{rgb}{0.25,0.88,0.82}
\definecolor{violetred}{rgb}{0.82,0.13,0.56}
\definecolor{violet}{rgb}{0.93,0.51,0.93}
\definecolor{wheat1}{rgb}{1.00,0.91,0.73}
\definecolor{wheat2}{rgb}{0.93,0.85,0.68}
\definecolor{wheat3}{rgb}{0.80,0.73,0.59}
\definecolor{wheat4}{rgb}{0.55,0.49,0.40}
\definecolor{wheat}{rgb}{0.96,0.87,0.70}
\definecolor{whitesmoke}{rgb}{0.96,0.96,0.96}
\definecolor{white}{rgb}{1.00,1.00,1.00}
\definecolor{yellow1}{rgb}{1.00,1.00,0.00}
\definecolor{yellow2}{rgb}{0.93,0.93,0.00}
\definecolor{yellow3}{rgb}{0.80,0.80,0.00}
\definecolor{yellow4}{rgb}{0.55,0.55,0.00}
\definecolor{yellowgreen}{rgb}{0.60,0.80,0.20}
\definecolor{yellow}{rgb}{1.00,1.00,0.00}

\DeclareMathAlphabet{\mathpzc}{OT1}{pzc}{m}{it} 

\theoremstyle{plain}
\newtheorem{thm}{Theorem}[section]
\newtheorem{cor}[thm]{Corollary}
\newtheorem{lmm}[thm]{Lemma}
\newtheorem{prpn}[thm]{Proposition}
\newtheorem{rem}[thm]{Remark}

\theoremstyle{definition}
\newtheorem{defn}[thm]{Definition}
\newtheorem{eg}[thm]{Example}

\newcommand*\Str{\vspace*{0.1cm}\includegraphics[scale=0.03]{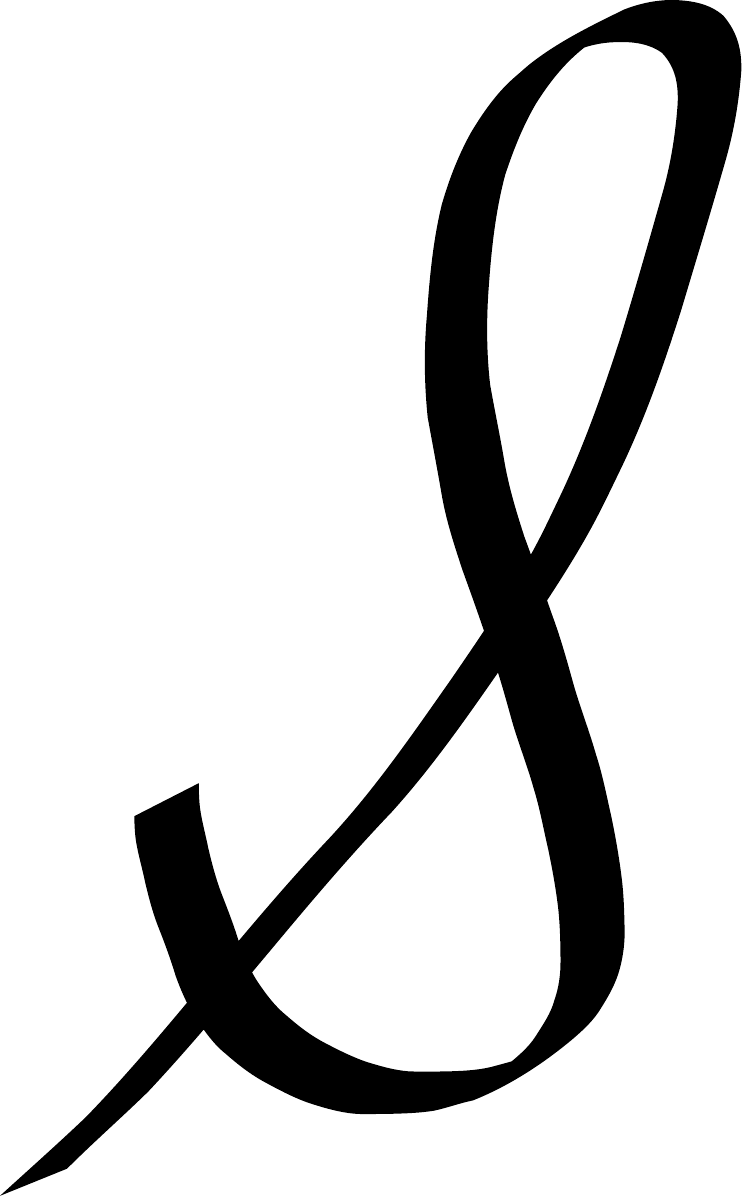}}

\numberwithin{equation}{section}

\newcommand{\N}{\mathbb{N}}

\newcommand{\Q}{\mathbb{Q}}
\newcommand{\R}{\mathbb{R}}

\newcommand{\Z}{\mathbb{Z}}

\newcommand{\into}{\hookrightarrow}

\newcommand{\lan}{\left\langle}
\newcommand{\ran}{\right\rangle}
\newcommand{\bgd}{\begin{displaymath}}
\newcommand{\edd}{\end{displaymath}}
\newcommand{\bge}{\begin{equation}}
\newcommand{\ede}{\end{equation}}
\newcommand{\bgea}{\begin{eqnarray*}}
\newcommand{\edea}{\end{eqnarray*}}
\newcommand{\bgeA}{\begin{eqnarray}}
\newcommand{\edeA}{\end{eqnarray}}
\newcommand{\bgc}{\begin{center}}
\newcommand{\edc}{\end{center}}
\newcommand{\ben}{\begin{enumerate}}
\newcommand{\een}{\end{enumerate}}
\newcommand{\bgi}{\begin{itemize}}
\newcommand{\edi}{\end{itemize}}
\newcommand{\hf}{\hspace*{0.5cm}}

\newcommand{\lra}{\longrightarrow}

\newcommand{\ep}{\varepsilon}
 
\newcommand{\lp}{\bullet}

\title{{\LARGE Realizing congruence subgroups inside the diffeomorphism group of a product of homotopy spheres}}
\author{Somnath Basu \and F. Thomas Farrell \thanks{supported in part by NSF grant DMS-1206622.}}
\date{}
\usepackage{pdfpages}
\usepackage[pdfauthor={Somnath Basu, F. Thomas Farrell},
        pdftitle={Realizing congruence subgroups inside diffeomorphism groups},
        pdftex, colorlinks, citecolor=red, filecolor=black, linkcolor=blue, urlcolor=black]{hyperref}

\begin{document}

\maketitle

\begin{abstract}
Let $M$ be a smooth manifold which is homeomorphic to the $n$-fold product of $S^k$, where $k$ is odd. There is an induced homomorphism from the group of diffeomorphisms of $M$ to the automorphism group of $H_k(M;\Z)$. We prove that the image of this homomorphism contains a congruence subgroup of $\textup{SL}_n(\Z)$ whenever $n$ is at least $3$.
\end{abstract}

\vspace*{0.5cm}
\tableofcontents
\vspace*{0.75cm}


\section{Introduction}

\hf\hf Given a group $G$ acting on a space, there are induced actions of $G$ on the endomorphisms of various topological invariants associated with the space. For instance, if $M$ is a manifold and $G$ acts on $M$ then there is an induced action of $G$ on the automorphisms of $H^\lp(M;\Z)$, the cohomology ring of $M$. A natural question is the following:
\bgc
{\it What is the image of $G$ inside $\textup{Aut}(H^\lp(M;\Z))$?}
\edc
 This question can be phrased in the appropriate category and the answer will depend on this choice. For example, one can choose $G$ to be acting by isometries on a Riemannian manifold. In this paper we will focus on groups acting by diffeomorphisms on a smooth manifold.\\

\hf\hf It is worthwhile mentioning Zimmer's program briefly (cf. \cite{Zim87, Zim87b} for Zimmer's original version of the conjectures). Although this program is not directly related to our paper, some of the results (of the program) that have been proven or have been conjectured resonate well with our result. Broadly construed, Zimmer's program aims to understand actions of finitely generated groups on compact manifolds by diffeomorphisms. A more precise version is the following conjecture: 
\bgd
\textit {Higher rank irreducible lattices do not act faithfully on low dimensional manifolds}.
\edd
In other words, Zimmer suggests that one should be able to show that $\textup{SL}_n(\Z)$ cannot act non-trivially on any manifold of dimension $d<n-1$ except via finite quotients. There is vast body of mathematical work and literature on Zimmer's program. We refer to the excellent review article \cite{Fis11} of Fisher for details and further references. In fact, a question posed in \cite{Fis11} (based on previous work by Fisher and Whyte) is the following:\\

{\bf Question}\hf {\it Let $M$ be a compact manifold with $\pi_1(M) \cong \Z^n$ and assume $\Gamma < \textup{SL}_n(\Z)$ has finite index. Let $\Gamma$ act on $M$ fixing a point so that the resulting $\Gamma$ action on $\pi_1 (M)$ is given by the standard representation of $\textup{SL}_n(\Z)$ on $\Z^n$. Is it true that the dimension of $M$ is at least $n$?}\\[0.2cm]

\hf\hf We formulate a variant of this question for simply connected manifolds. This variant deviates, in part, from Zimmer's program in the sense that we are trying to construct lattices that act faithfully via diffeomorphisms on manifolds whose dimension is much higher than the rank of the lattice. One may consider a smooth manifold $M$ which is {\it homeomorphic} to the product of $n$ smooth spheres, i.e.,
\bgd
f:M\stackrel{\simeq}{\lra} \Sigma_1^k\times\Sigma_2^k\times\cdots\times \Sigma_n^k
\edd
is a homeomorphism and $\Sigma_j^k$ are smooth homotopy $k$-spheres. (By a homotopy $k$-sphere we mean any smooth manifold which is homeomorphic to $S^k$.) The group of diffeomorphisms of $M$ induces an action on $H^k(M;\Z)\cong \Z^n$. In fact, this action determines the action on the full cohomology ring by the cup product structure. We shall focus on actions on homology groups of $M$. One can ask if there is a finite index subgroup $\Gamma < \textup{SL}_n(\Z)$ which acts by diffeomorphisms on $M$ and induces the standard representation of $\Gamma<\textup{SL}_n(\Z)$ on $\Z^n\cong H_k(M;\Z)$? It is not apriori clear that such actions exist. The following is our main result.\\
%

{\bf Theorem}\hf {\it Let $M$ be a smooth manifold homeomorphic to a $n$-fold product of $k$-dimensional smooth spheres where $k$ is odd. The subgroup $\mathfrak{R}(k,n)$ of $\textup{GL}_n(\Z)$ representable by self-diffeomorphisms of $M$ contains a congruence subgroup if $k=1,n\neq 4$ or $k>1,n\geq 3$.}\\[0.4cm]
When $M$ is a product of standard spheres then a similar result was proven by Lucas and Saeki \cite{LuSa02}. However, our result is far more general and the method of proof is necessarily different. \\

\hf\hf Several comments are in order. Our main result is a combination of Theorem \ref{mainthm}, Theorem \ref{mainthm2} and Remark \ref{mainthm3}. Theorem \ref{mainthm} is a special case when $M$ {\it is} the product of $n$ smooth $k$-spheres and $k>1$. A modified line of reasoning leads us to Remark \ref{mainthm3} which allows us to upgrade Theorem \ref{mainthm} by replacing a product of spheres up to homeomorphism. The hypothesis on $M$ in Theorem \ref{mainthm2} (corresponding to $k=1$ in Theorem) is that $M$ is homotopy equivalent to the $n$-torus. It is known that a manifold homotopy equivalent to the $n$-torus is actually homeomorphic to it. \\

\hf\hf Congruence subgroups, being the kernel of the natural map\footnote{The map reduces the coefficients modulo $m$; here $\Z_m:=\Z/m\Z$.} $\textup{SL}_n(\Z)\to \textup{SL}_n(\Z_m), m\in\N$, are of finite index. Therefore, our theorem provides a step towards an affirmative answer to the question posed earlier. However, ideally one would want to realize this congruence subgroup inside the diffeomorphism group; i.e., split $\mathfrak{R}(k,n)$, or some congruence subgroup inside it, back to $\textup{Diff}(M)$. And this is a natural direction we hope to pursue.\\[0.2cm]

\section{Monoid of homotopy equivalences}

\hf\hf Let $S^k$ be a simply connected sphere, where the $k$ is odd and $k\neq 1,3,7$. The cases of $S^1, S^3$ and $S^7$ will be dealt with in \S \ref{osph}. Let
\bgd
f:M\stackrel{\simeq}{\lra}M(k,n):=\Sigma_1^k\times\cdots\times \Sigma_n^k
\edd
be a smooth manifold homeomorphic to the product of $n$ smooth spheres of dimension $k$. We shall use $\mathbb{S}^k$ to denote the {\it round sphere}, i.e., the unit sphere in $\R^{k+1}$ equipped with the standard smooth structure. 
\begin{defn}
Given an oriented, topological manifold $X$, we define $\textup{Aut}(X)$ to be the set of orientation preserving homotopy equivalences from $X$ to itself. The set of homotopy classes of (orientation preserving) self-homotopy equivalences of $X$ forms a group and will be denoted by $\mathcal{G}(X)$.
\end{defn}
\begin{rem}
Since $\Sigma_i^k$ and $\mathbb{S}^k$ have the same underlying topological manifold, for the purpose of studying $\textup{Aut}(M)$, it suffices to consider self-homotopy equivalences of $(S^k)^n=S^k\times\cdots\times S^k$.
\end{rem}
\hf\hf Let $\iota_j:S^k\into (S^k)^n$ be the inclusion of the $j^\textup{th}$ sphere. We shall also use $\iota_j$ to denote the image of a generator $\iota_j[S^k]\in\pi_k\big((S^k)^n\big)$. It follows from the definition of $(S^k)^n$ that the Whitehead products $[\iota_j,\iota_l]=0$ whenever $j\neq l$. We shall assume that $n\geq 2$. Given a map $f:(S^k)^n\to (S^k)^n$, the induced map
\bge\label{ind}
f_\ast :H_k\big((S^k)^n;\Z\big)\lra H_k\big((S^k)^n;\Z\big),
\ede
in terms of the standard basis coming from the product decomposition, is an element of $\textup{M}_n(\Z)$. We shall denote this $n\times n$ integer matrix by $A(f)=((a_{ij}))$. Let $\iota\in\pi_k(S^k)$ be the class of the identity map. It is known \cite{Ada60} that $[\iota,\iota]$ is of order $2$ when $k$ is odd and not equal to $3$ or $7$. Then
\begin{eqnarray*}
0=f_\ast[\iota_j,\iota_l] & = & (a_{j1}a_{l2}-a_{j2}a_{l1})[\iota_1,\iota_2]+\cdots+(a_{j(n-1)}a_{ln}-a_{jn}a_{l(n-1)})[\iota_{n-1},\iota_n]\\
& & + a_{j1}a_{l1}[\iota_1,\iota_1]+\cdots+ a_{jn}a_{ln}[\iota_n,\iota_n]  \\
& = & \big(a_{j1}a_{l1}[\iota_1,\iota_1],\cdots,a_{jn}a_{ln}[\iota_n,\iota_n]\big)\in \Z_2^{\oplus n}<\pi_{2k-1}\big((S^k)^n\big).
\end{eqnarray*}
This forces $a_{js}a_{ls}$ to be even for $j\neq l$ and any $s$. 
\begin{defn}
Given two elements $\mathbf{v},\mathbf{w}\in\Z^n$ we define the {\it pre-dot product} to be
\bgd
\mathbf{v}\circ\mathbf{w}:=(v_1w_1,\cdots,v_nw_n)\in\Z^n.
\edd
Let $\mathpzc{W}_n(2)$ denote the set of matrices $A\in\textup{SL}_n(\Z)$ with the property that the pre-dot product of any two distinct rows of $A$ lie in $(2\Z)^n\subset\Z^n$.
\end{defn}
\begin{lmm}
The set $\mathpzc{W}_n(2)$ is a subgroup of $\textup{SL}_n(\Z)$ containing the principal congruence subgroup $\Gamma_n(2)$ with $[\mathpzc{W}_n(2):\Gamma_n(2)]=n!$. A set of coset representatives is given by $A_n\sqcup \tau A_n$ with $\tau$ as in \eqref{tau} and $A_n$ is the group of even permutations.
\end{lmm}
{\bf Proof.}\hf It is known that the homomorphism (reduction of coefficients)
\bgd
R:\textup{SL}_n(\Z)\lra \textup{SL}_n(\Z_2)
\edd
is surjective. The set $P_n$ of permutations on $n$ letters is a subgroup\footnote{All permutations have sign $1$ in $\Z_2$.} of $\textup{SL}_n(\Z_2)$. We shall show that $\mathpzc{W}_n(2)=R^{-1}(P_n)$. In particular, $\Gamma_n(2):=R^{-1}(\textup{Id})$ is a normal subgroup of $\mathpzc{W}_n(2)$. We are led to a short exact sequence
\bgd
1\lra \Gamma_n(2)\lra \mathpzc{W}_n(2)\stackrel{R}{\lra} P_n\lra 1.
\edd
Therefore, the index of $\Gamma_n(2)$ in $\mathpzc{W}_n(2)$ is $n!$. \\
\hf It is clear that $\mathpzc{W}_n(2)\subseteq R^{-1}(P_n)$. Now let $A\in \mathpzc{W}_n(2)$. Since $\det A$ is odd, each column of $A$ must have at least one odd entry, for otherwise $\det A$ (via a cofactor expansion along this column) is even. By the defining property of $A$, each column of $A$ can have at most one odd entry, whence exactly one such entry per column is permitted. Thus, $A$ has exactly $n$ odd entries. If a row has more than one odd entry then some row has all entries even which contradicts $\det A$ being odd; $A$ must have exactly one odd entry in each row and each column, whence $R(A)\in P_n$.\\
\hf Let $P_\sigma$ be the permutation matrix corresponding to a permutation $\sigma$ of $\{1,\cdots,n\}$, i.e., $(P_\sigma)_{ij}=\delta_{j\sigma(i)}$. Let $\tau\in \mathpzc{W}_n(2)$ be the matrix
\bge\label{tau}
\tau:=\left(\begin{array}{cc|c}
0 & -1 &  \\
1 & 0 &  \\
\hline 
 & & I_{n-2}
\end{array}\right).
\ede 
This plays the role of the transposition\footnote{We cannot use the usual permutation matrix associated to a transposition as it has determinant $-1$.} $(1\,2)$. It can be verified that 
\bgd
\mathpzc{W}_n(2)=\big(\sqcup_{\sigma\in A_n}P_\sigma\Gamma_n(2)\big)\sqcup\big(\sqcup_{\sigma\in A_n}\tau P_\sigma\Gamma_n(2)\big).
\edd
In particular, this implies that $A_n$ and $\tau A_n$ form a set of coset representatives for $\Gamma_n(2)$ in $\mathpzc{W}_n(2)$.$\hfill\square$\\[0.2cm]
\hf\hf We now consider orientation preserving homotopy equivalences of $M(k,n)$. Let $f$ be a homotopy equivalence with $A(f)$ being the induced matrix \eqref{ind}. Since $f$ is orientation preserving, the induced matrix $A(f)$ has determinant $1$. Consider the map
\bgd
\Psi:\textup{Aut}\big(M(k,n)\big)\lra \mathpzc{W}_n(2),\,\,f\mapsto A(f).
\edd
This map induces a group homomorphism 
\bgd
\Psi:\mathcal{G}\big(M(k,n)\big)\lra \mathpzc{W}_n(2),
\edd
where $\mathcal{G}(X)$ is the group of (orientation preserving) homotopy classes of self-homotopy equivalences.
\begin{prpn}\label{wn}
The group homomorphism $\Psi:\mathcal{G}\big(M(k,n)\big)\lra \mathpzc{W}_n(2)$ is surjective.
\end{prpn}
\begin{rem}
Any element $A\in \mathpzc{W}_2(2)$ can be realized by a map; define $f_A:S^k\vee S^k\to S^k\vee S^k$ by wrapping the first sphere $a_{11}$ times around the first sphere and $a_{12}$ times around the second sphere. Now wrap the second sphere $a_{21}$ times around the first sphere and $a_{22}$ times around the second sphere. Since $a_{11}a_{21}$ and $a_{21}a_{22}$ are both even, 
\bgd
f_A[\iota_1,\iota_2]=(a_{11}a_{22}-a_{21}a_{21})[\iota_1,\iota_2]+a_{11}a_{21}[\iota_1,\iota_1]+a_{12}a_{22}[\iota_2,\iota_2]=[\iota_1,\iota_2].
\edd
We extend $f_A$ to $f:S^k\times S^k\to S^k\times S^k$ by attaching a top cell on both sides, where the attaching map is the same on both sides. By construction, $f\in\textup{Aut}\big(M(k,2)\big)$ and induces $A(f)=A$.
\end{rem}
{\bf Proof.}\hf To prove surjectivity of $\Psi$, it suffices to check surjectivity for a generating set of $\mathpzc{W}_n(2)$. For $n\geq 3$ each even permutation $\sigma\in A_n$ naturally induces (and is induced by) the map that permutes the factors in $M(k,n)$. To realize $\tau$, we can define a map $f:S^k\vee S^k\lra S^k\vee S^k$ that induces the non-trivial $2\times 2$ block on $H_k(S^k\vee S^k;\Z)$. Note that the Whitehead product $[\iota_1,\iota_2]$ is mapped to $[\iota_1,\iota_2]$ under $f$. We can attach a cell of dimension $2k$ on both sides along $[\iota_1,\iota_2]$ extending $f$ to an automorphism of $S^k\times S^k$.  Setting 
\bgd
f_\tau:=f\times \textup{id}^{(n-2)}:M(k,n)\lra M(k,n)
\edd
implies that $A(f_\tau)=\tau$. This takes care of the generators of the coset representatives.\\
\hf Recall that we use $\ep_{ij}$ to denote the matrix with zeroes everywhere except at $(i,j)^{\textup{th}}$ entry which is $1$, whence the elementary matrix $E_{ij}=I_n+\ep_{ij}$. Lemma \ref{gen} implies that $\Gamma_n(2)$ is generated by:\\
\hf (i) the second power $E_{ij}^2$ (and their inverses) of elementary matrices $E_{ij}$, and\\
\hf (ii) the matrices $J_i=I_n-2\ep_{ii}-2\ep_{(i+1)(i+1)}$.\\
For any $E_{ij}^2$ the matrix $(1,2;0,1)$ or $(1,0;2,1)$ can be realized by a self-homotopy equivalence of $S^k\times S^k$ (the product of the $i^\textup{th}$ and $j^\textup{th}$ sphere in $M(k,n)$). This can be extended (by declaring this to be identity on the other factors) to an element of $f\in\textup{Aut}\big(M(k,n)\big)$ such that $A(f)=E_{ij}^2$. The non-trivial $2\times 2$ block in $J_i$ is the matrix $-I_2\in \mathpzc{W}_2(2)$. This can be realized by a self-homotopy equivalence of $S^k\times S^k$ which changes orientation of $i^\textup{th}$ and $(i+1)^\textup{th}$ spheres simultaneously. This map extends (via identity on other factors) to an element of $\textup{Aut}\big(M(k,n)\big)$. This proves surjectivity of $\Psi$.  $\hfill\square$\\[0.2cm]
It still remains to prove the following result.
\begin{lmm}\label{gen}
\textup{(i)} The principal congruence subgroup $\Gamma_2(2)$ is generated by $E_{12}^2, E_{21}^2$ and $-I_2$. \\
\textup{(ii)} When $n\geq 3$ the principal conguence subgroup $\Gamma_n(2)\triangleleft\,\textup{SL}_n(\Z)$ is generated by $E_{ij}^2$ and 
\bgd\label{diag}
J_i=\left(\begin{array}{c|rr|c}
I_{i-1} & & & \\
\hline
& -1 & 0 &  \\
& 0 & -1 &  \\
\hline
& & & I_{n-i+1} \\
\end{array}\right)
\edd
for $1\leq 1< n$.
\end{lmm}
{\bf Proof.}\hf We shall prove (ii) first, which is a generalization of (i). Let 
\bgd
S:=\{E_{ij}^{2}\,|\,1\leq i,j\leq n,i\neq j\}\subset \Gamma_n(2). 
\edd
It is known \cite{BMS67} that the normal closure of the subgroup generated by $S$ in $\textup{SL}_n(\Z)$ is $\Gamma_n(2)$ when $n\geq 3$. However, $S$ itself does not generate $\Gamma_n(2)$.
Let $H$ be the subgroup generated by $S$ and $J_i$'s. The normal closure of $H$ in $\textup{SL}_n(\Z)$ is necessarily $\Gamma_n(2)$. The proof is complete if we show that $H$ is normal. Equivalently, since $E_{ij}$'s generate $\textup{SL}_n(\Z)$, we need to show that $E_{ij}AE_{ij}^{-1}$ and $E_{ij}^{-1}A E_{ij}$ are words in $E_{ij}^2$'s and $J_i$'s when $A\in S\cup \{J_i\}_i$. If we denote the matrix $J_{ik}=I_n-2\ep_{ii}-\ep_{kk}$ then
\bgd
J_{ik}=\left\{\begin{array}{cc}
J_i J_{i+1}\cdots J_k & \textup{if $i<k$}\\
J_k J_{k+1}\cdots J_i  & \textup{if $i>k$}
\end{array}\right.
\edd
is an element of $H$. In this notation, $J_{i(i+1)}=J_i$ as defined before. A few explicit calculations reveal that
\bgd
E_{ij} E_{kl}^2 E_{ij}^{-1}=\left\{\begin{array}{rl}
E_{kl}^2 & \textup{if $j\neq k,i\neq l$}\\
E_{kl}^2E_{kj}^{-2}  & \textup{if $j\neq k,i=l$}\\
E_{il}^2E_{kl}^{2}  & \textup{if $j=k,i\neq l$}\\
E_{ik}^2E_{ki}^{-2}J_{ik}   & \textup{if $j=k,i=l$} \\
\end{array}\right.
\edd
\bgd
E_{ij}^{-1} E_{kl}^2 E_{ij}=\left\{\begin{array}{rl}
E_{kl}^2 & \textup{if $j\neq k,i\neq l$}\\
E_{kl}^2E_{kj}^{2}  & \textup{if $j\neq k,i=l$}\\
E_{il}^{-2}E_{kl}^{2}  & \textup{if $j=k,i\neq l$}\\
J_{ki}E_{ki}^{-2}E_{ik}^2  & \textup{if $j=k,i=l$} \\
\end{array}\right.
\edd
\bgd
E_{ij} J_{k} E_{ij}^{-1}=\left\{\begin{array}{rl}
J_{k} & \textup{if $\{i,j\}\cap \{k,k+1\}=\varnothing$}\\
J_k E_{ij}^{-2}  & \textup{if $i\in \{k,k+1\}, j\not\in\{k,k+1\}$}\\
E_{ij}^{2}J_k  & \textup{if $i\not\in \{k,k+1\}, j\in\{k,k+1\}$}\\
J_{k} & \textup{if $\{i,j\}\cap \{k,k+1\}=\{k,k+1\}$}
\end{array}\right.
\edd
\bgd
E_{ij}^{-1} J_{k} E_{ij}=\left\{\begin{array}{rl}
J_{k} & \textup{if $\{i,j\}\cap \{k,k+1\}=\varnothing$}\\
E_{ij}^{-2}J_k  & \textup{if $i\in \{k,k+1\}, j\not\in\{k,k+1\}$}\\
J_kE_{ij}^{2}  & \textup{if $i\not\in \{k,k+1\}, j\in\{k,k+1\}$}\\
J_{k} & \textup{if $\{i,j\}\cap \{k,k+1\}=\{k,k+1\}$}.
\end{array}\right.
\edd
\hf To prove (i), which is well-known and has several proofs, we recall that any upper triangular matrix with $\pm 1$'s on the diagonal is generated by $E_{12}^2$ or $E_{21}^2$. For any other matrix $A=(a,2b;2c,d)$ in $\Gamma_2(2)$ we may use the Euclidean algorithm to reduce the modulus of either $|a|$ or $|d|$ by multiplying on the right by a suitable power of $E_{12}^2$ or $E_{21}^2$; this is possible since $\max\{|a|,|d|\}>\min\{2|b|,2|c|\}$. An iteration of this leads to one of the diagonal entry being $1$; such matrices are of the form $E_{12}^{2b} E_{21}^{2c}$ or $E_{21}^{2c}E_{12}^{2b}$.       $\hfill\square$
\begin{defn}
Let $h\mathfrak{R}(k,n)$ denote the subgroup of $\textup{GL}_n(\Z)$ representable by a self-homotopy equivalence $f$ of $M(k,n)$, i.e., $A\in h\mathfrak{R}(k,n)$ if and only if 
\bgd
A(f)=f_\ast:H_k\big(M(k,n);\Z\big)\to H_k\big(M(k,n);\Z\big),
\edd
in terms of the natural basis, is $A$. 
\end{defn}
In light of Proposition \ref{wn}, we have a description of $h\mathfrak{R}(k,n)$.
\begin{thm}\label{hrnk}
The group $h\mathfrak{R}(k,n)$ contains $\mathpzc{W}_n(2)$ as a subgroup of index $2$.
\end{thm}
{\bf Proof.}\hf It follows from Proposition \ref{wn} that $\mathpzc{W}_n(2)$ is precisely the subgroup of $h\mathfrak{R}(k,n)$ consisting of orientation preserving homotopy equivalences. This is subgroup of index $2$. In fact, a precise description of $h\mathfrak{R}(k,n)$ is that it consists of elements $A\in \textup{GL}_n(\Z)$ such that $A\!\mod 2$ is a permutation matrix.$\hfill\square$
\begin{rem}\label{stndsp}
It is a classical fact that the map $\psi_\mathbf{x}:S^k \to S^k$, defined for $\mathbf{x}\in S^k$, as
\bgd
\psi_\mathbf{x}(\mathbf{y}):=\mathbf{x}-2\lan \mathbf{x},\mathbf{y}\ran \mathbf{y}
\edd
has degree $(-1)^{k+1}+1$. If $k$ is odd then this degree is $2$. In \cite{LuSa02} they use this to realize the matrix $E_{12}^2$ via the diffeomorphism
\bgd
f:(\mathbb{S}^k)^n\lra (\mathbb{S}^k)^n,\,\,(\mathbf{x}_1,\mathbf{x}_2,\cdots,\mathbf{x}_n)\mapsto \big(\psi_{\mathbf{x}_1}(\mathbf{x}_2),\mathbf{x}_2,\cdots, \mathbf{x}_n\big).  
\edd
The matrices $E_{ij}^2$ can be realized in a similar manner. Therefore, any element of $h\mathfrak{R}_\textup{st}(k,n)$ can actually be realized by diffeomorphisms. 
\end{rem}

\subsection{The other spheres}\label{osph}

\hf\hf It remains to discuss the cases of self-homotopy equivalences of $M(k,n)$ when either $k$ is even or $k=1,3,7$. When $k$ is even, we shall see that the group $h\mathfrak{R}(k,n)$ is actually finite. The cases when $k=1,3,7$ are special from other odd values of $k$ as $S^1, S^3$ and $S^7$ admit multiplications. \\
\hf\hf We first deal with the case when $M(k,n)=\Sigma_1^k\times\cdots\times \Sigma_n^k$ and $k$ is even. Note that if $\iota\in \pi_k(S^k)$ then $[\iota,\iota]\in\pi_{2k-1}(S^k)$ is of infinite order. Let $f:M(k,n)\to M(k,n)$ be a map with $f_\ast=A(f)=((a_{ij}))$ on $H_k\big(M(k,n);\Z\big)$. Let $\iota_j:S^k\into M(k,n)$ be the inclusion of the $j^\textup{th}$ sphere. We shall also use $\iota_j$ to denote the homology class $(\iota_j)_\ast[S^k]\in H_k\big(M(k,n);\Z\big)$.  Since $[\iota_j,\iota_l]=0$ when $j\neq l$, we have
\bgd
0=f_\ast[\iota_j,\iota_l]=[f_\ast(\iota_j),f_\ast(\iota_l)]=\sum_{k=1}^na_{jk}a_{lk}[\iota_k,\iota_k].
\edd
In particular, each of the coefficients $a_{jk}a_{lk}=0$ and in the $k^\textup{th}$ column of $A$ there is at most one non-zero entry. Thus, each column of $A$ can only have at most one non-zero entry. If we demand that $f$ is a homotopy equivalence then $\det A=\pm 1$ and $A$, being a matrix with integer entries, has only $\pm 1$ as its non-zero entries. Such matrices are realizable as maps from $M(k,n)$ to itself which permutes the factors and (possibly) switch orientation of some factors. The group $h\mathfrak{R}(k,n)\subset \textup{GL}_n(\Z)$ is, in fact, of order $2^n n!$.\\
\hf\hf Now we turn to $M(k,n)$ with $k=1,3,7$. The spheres $S^1,S^3$ and $S^7$ admit multiplication. In fact, $\mathbb{S}^1, \mathbb{S}^3$ are Lie groups but $S^7$ does not admit a strictly associative multiplication \cite{Jam57}. Let $\iota\in \pi_k(S^k)$ be the class of the identity map. Then $[\iota,\iota]=0$ and this happens \cite{Ada60} only for $k=1,3,7$. Therefore, there are no obstructions from Whitehead products as far as realizability of a matrix $A\in \textup{GL}_n(\Z)$ is concerned. In fact, given $A\in \textup{GL}_n(\Z)$ and $k=1,3$, let us define
\bge\label{PA}
\mathpzc{P}_A:(\mathbb{S}^k)^n\lra (\mathbb{S}^k)^n,\,\,(x_1,\cdots, x_n)\mapsto (x_1^{a_{11}}x_2^{a_{12}}\cdots x_n^{a_{1n}},\cdots,x_1^{a_{n1}}x_2^{a_{n2}}\cdots x_n^{a_{nn}}),
\ede
where we consider $\mathbb{S}^3$ as the multiplicative group $\textup{SU}(2)$. This is actually a real analytic map and induces the matrix $A$. Thus, there is the equality
\bge\label{hRkn13}
h\mathfrak{R}(k,n)=\textup{GL}_n(\Z),\,\,\,k=1,3.
\ede
Note that the definition of $\mathpzc{P}_A$ does not work for $S^7$ due to the failure of associativity \cite{Jam57} of any multiplication on $S^7$.
\begin{rem}
Notice that $\mathpzc{P}_A$ \eqref{PA} is a diffeomorphism for $(\mathbb{S}^1)^n$. However, the map $\mathpzc{P}_A$ fails to be a diffeomorphism  for $(\mathbb{S}^3)^n$ for most matrices $A\in\textup{GL}_n(\Z)$. An explicit example where $\mathpzc{P}_A$ is not one-to-one is the $2\times 2$ matrix $A$ whose first row is $(1,-1)$ and second row is $(-1,2)$. In this case, both $(-i,-1)$ and $\textstyle{\frac{1}{2}}(i+ak,1+aj)$ are mapped to
$(i,i)$ where $a$ denotes a square root of $3$ and $i,j,k$ are the standard unit quaternions.
\end{rem}
\hf\hf The group $\textup{SL}_n(\Z)$ is generated by the elementary matrices $E_{ij}$. Now any $E_{ij}$ is realizable for $(\mathbb{S}^1)^n, (\mathbb{S}^3)^n$ and $(\mathbb{S}^7)^n$, where $\mathbb{S}^7$ is the space of unit octonions. For instance, if $i<j$ then the map 
\bgd
\mathpzc{P}_{ij}(x_1,\cdots,x_n)=(x_1,\cdots, x_{i-1},x_ix_j,x_{i+1},\cdots,x_n)
\edd
is one such, i.e., $A(\mathpzc{P}_{ij})=E_{ij}$. It is clear that $\mathpzc{P}_{ij}$ is a diffeomorphism. Given $A\in\textup{SL}_n(\Z)$ we may write it as a product of elementary matrices $A=E_1E_2\cdots E_r$. Then $\mathpzc{P}_{E_1}\circ\cdots\circ \mathpzc{P}_{E_r}$ is a diffeomorphism which induces $A$. In general, since exotic spheres exist in dimension $7$, it is not clear whether $E_{ij}$ can be realized for $M(7,n)=\Sigma_1^7\times\cdots\times\Sigma_n^7$. For the convenience of the reader, we gather all the previous observations.
\begin{prpn}\label{othsph}
Let $M(k,n)$ be as before and $h\mathfrak{R}(k,n)$ be the associated group.\\
\textup{(1)} If $k$ is even then $h\mathfrak{R}(k,n)$ is a subgroup of $\textup{GL}_n(\Z)$ of order $2^n n!$.\\
\textup{(2)} If $k=1,3,7$ then $h\mathfrak{R}(k,n)=\textup{GL}_n(\Z)$.
\end{prpn}

\section{Subgroups in the image of diffeomorphism groups}

\hf\hf Recall our notation: let $k\geq 3$ be an odd positive integer and $n\geq 1$ be a positive integer, and let 
\bgd
M(k,n)=\Sigma_1^k\times \cdots\times \Sigma_n^k,
\edd
where $\Sigma_i^k$ is a {\it smooth} sphere of dimension $k$.
\begin{lmm}\label{ssfin}
The set $\textup{\Str}_0\big(M(k,n)\big)$ is finite. 
\end{lmm}
Here $\textup{\Str\,}\big(M(k,n)\big)$ is the {\it smooth structure set} of $M^{kn}:=M(k,n)$, i.e., each element in $\textup{\Str\,}(M^{kn})$ is represented by a homotopy equivalence $f:N^{kn}\to M^{kn}$, where $N^{kn}$ is a closed, smooth manifold - called a {\it homotopy smooth structure} on $M^{kn}$. If $g:K^{kn}\to M^{kn}$ is a second such structure, then $f\equiv g$ if and only if there exists a diffeomorphism $F:N^{kn}\to K^{kn}$ such that $g\circ F\sim f$, i.e., $g\circ F$ is homotopic to $f$. Note that $\textup{\Str\,}(M^{kn})$ has an obvious base point $\ast$, namely, $\ast$ is the class of $\textup{id}_{M}:M^{kn}\to M^{kn}$. One defines $\textup{\Str}_0(M^{kn})$ to be the subset of $\textup{\Str\,}(M^{kn})$ consisting of those elements that can be represented by some self homotopy equivalence $\varphi:M^{kn}\to M^{kn}$. It is clear that $\ast\in\textup{\Str}_0(M^{kn})$. If $M$ and $N$ are smooth manifolds which are homotopy equivalent then there is a bijection between $\textup{\Str}(M)$ and $\textup{\Str}(N)$.
\begin{eg}
For $l\geq 1$ we note that $\textup{\Str\,}(\mathbb{S}^l)$ is actually the group of exotic spheres $\Theta_l$ and $\textup{\Str}_0(\mathbb{S}^l)=\{\ast\}.$ Moreover, the structure sets $\textup{\Str\,}(\mathbb{S}^3\times \mathbb{S}^4)$ and $\textup{\Str\,}(\mathbb{S}^3\times \mathbb{S}^3\times \mathbb{S}^3\times \mathbb{S}^3\times \mathbb{S}^3)$ are both infinite.
\end{eg}
\begin{cor}\label{Rnk}
Let $m:=m(k,n)=\big|\textup{\Str}_0(M^{kn})\big|$ and $\varphi:M\to M$ be a self-homotopy equivalence. Then there exists a self-diffeomorphism $f:M\to M$ such that $f\sim \varphi^{m!}$.
\end{cor}
{\bf Proof.}\hf Consider the set $\{\textup{id}=\varphi^0,\varphi,\cdots, \varphi^m\}$ of cardinality $m+1$. By the pigeon-hole principle, there exists integers $0\leq i<j=i+s\leq m$ such that the homotopy smoothing $\varphi^i:M\to M$ and $\varphi^j:M\to M$ are equivalent, i.e., there exists a self-diffeomorphism $F:M\to M$ such that the following diagram commutes up to homotopy:
\bgd
\xymatrix{
M\ar[r]^-{\varphi^i} & M\\
M\ar[u]^-{F}\ar[ur]_-{\varphi^{j}} & }
\edd
where $j=i+s$. Since $(\varphi^i\circ\varphi^s)\sim \varphi^i\circ F$ we conclude that $F\sim \varphi^s$, whence $F^{\frac{m!}{s}}\sim \varphi^{m!}$. $\hfill\square$
\begin{defn}
Let $\mathfrak{R}(k,n)$ denote the subgroup of $\textup{GL}_n(\Z)$ representable by $\textup{Diff}\big(M(k,n)\big)$, i.e., $A\in \mathfrak{R}(k,n)$ if and only if there exists a diffeomorphism $f:M(k,n)\to M(k,n)$ such that 
\bgd
f_\ast:H_k\big(M(k,n);\Z\big)\lra H_k\big(M(k,n);\Z\big)
\edd
is given\footnote{Recall that we are using the natural integral basis of $H_k\big(M(k,n);\Z\big)$ arising from the product decomposition.} by the matrix $A$.
\end{defn}
\begin{rem}
Since $h\mathfrak{R}(k,n)$ contains $\mathfrak{R}(k,n)$, it follows from Proposition \ref{othsph} and Theorem \ref{hrnk} that in general all of $\textup{GL}_n(\Z)$ is not representable by diffeomorphisms. In fact this is only possible when  $k = 1,3$ or $7$.
\end{rem}

We shall reserve the notation $\mathfrak{R}_{\textup{st}}(k,n)$ for the group $\mathfrak{R}(k,n)$ associated to $\mathbb{S}^k\times\cdots\times \mathbb{S}^k$, where $\mathbb{S}^k$ is the standard smooth sphere.
\begin{cor}\label{Rknnor}
The group $\mathfrak{R}(k,n)$ contains an infinite normal subgroup of $\textup{GL}_n(\Z)$.
\end{cor}
To prove this we make crucial use of the related larger subgroup $h\mathfrak{R}(k,n)$ consisting of matrices $A\in \textup{GL}_n(\Z)$ such that there exists a self-homotopy equivalence $f$ of $M(k,n)$ with $f_\ast=A$. We will relate $\mathfrak{R}(k,n)$ to the congruence subgroup $\Gamma_n(2)$ of $\textup{GL}_n(\Z)$. Recall that for $s\in\Z^+$, the congruence subgroup $\Gamma_n(s)$ of level $s$ is the kernel of the natural homomorphism $\textup{SL}_n(\Z)\to\textup{SL}_n(\Z_s)$ induced by the ring homomorphism $\Z\to \Z_s$. Using the fact that $\textup{SL}_n(\Z_s)$ is a finite group, it is easily seen that $\Gamma_n(s)$ is a normal subgroup having finite index in $\textup{GL}_n(\Z)$. \\
\hf\hf To prove Corollary \ref{Rknnor}, we recall our previously proven result (cf. Theorem \ref{hrnk}).
\begin{thm}
The group $h\mathfrak{R}(k,n)$ contains $\Gamma_n(2)$.
\end{thm}
{\bf Proof of Corollary \ref{Rknnor}.}\hf Let $S=\{g^{m!}\,|\,g\in\Gamma_n(2)\}$, where $m=\big|\textup{\Str\,}_0(M^{kn})\big|$. We will show that the subgroup $\overline{S}$, generated by $S$, is an infinite subgroup of $\mathfrak{R}(k,n)$ which is also a normal subgroup of $\textup{GL}_n(\Z)$; thus proving Corollary \ref{Rknnor}. The group $\overline{S}$ is clearly infinite (since it contains the elementary matrix $E_{12}^{2m!}$ and its powers). It is a normal subgroup of $\textup{GL}_n(\Z)$ because of the following elementary group theory result.\\[0.2cm]
{\bf Fact.}\hf {\it Let $N$ be a normal subgroup of $G$ and $t\in\Z^+$. The subgroup $N^t$ of $N$ generated by $S=\{a^t\,|\,a\in N\}$ is a normal subgroup of $G$.}\\[0.2cm]
{\bf Proof of fact.}\hf Let $g\in G$ and $a^t\in S$, where $a\in N$. As $N$ is normal in $G$, $gag^{-1}\in N$ and, therefore, $ga^t g^{-1}=(gag^{-1})^t\in S$.\\[0.2cm]
Finally, $S\subseteq \mathfrak{R}(k,n)$ because of Corollary \ref{Rnk}. This completes the proof of Corollary \ref{Rknnor}. $\hfill\square$\\[0.2cm]
\hf\hf Now Margulis \cite{mar91} showed that every infinite normal subgroup of $\textup{SL}_n(\Z)$ (for $n\geq 3$) has finite index. Also Bass, Lazard and Serre \cite{BLS64} and Mennicke \cite{Men65} proved that every subgroup of finite index in $\textup{SL}_n(\Z)$, when $n\geq 3$, contains a congruence subgroup. Combining these deep results with Corollary \ref{Rknnor}, we obtain our second theorem.
\begin{thm}\label{mainthm}
The subgroup $\mathfrak{R}(k,n)$ of $\textup{GL}_n(\Z)$ representable by self-diffeomorphisms of $\,\Sigma_1^k\times\cdots\times \Sigma_n^k$ contains a congruence subgroup whenever $k$ is odd and $n\geq 3$.
\end{thm}
A similar line of argument leads to our next result.
\begin{thm}\label{mainthm2}
Let $\mathcal{T}^n$ be a homotopy $n$-torus, i.e., a smooth manifold which is homotopy equivalent to $\mathbb{T}^n$. The subgroup $\mathfrak{R}(\mathcal{T}^n)$ of $\textup{GL}_n(\Z)$ representable by self-diffeomorphisms of $\mathcal{T}^n$ contains a congruence subgroup whenever $n\neq 4$.
\end{thm}
{\bf Proof.}\hf When $n=1,2$ the homotopy torus $\mathcal{T}^n$ is diffeomorphic to $\mathbb{T}^n$ and the group $\mathfrak{R}(\mathcal{T}^n)$ is $\textup{GL}_n(\Z)$. When $n=3$, the proof of Poincar\'{e} conjecture by Perelman implies that homotopy $3$-tori are irreducible, i.e., any embedded $2$-sphere bounds a $3$-ball. Stalling's fibering Theorem \cite{Sta62} then implies that $\mathcal{T}^3$ is the total space of fibration with base $S^1$ and fibre $\mathbb{T}^2$. It follows from \cite{Neu63} that $\mathcal{T}^3$ is homeomorphic to $\mathbb{T}^3$ and by Moise's theorem \cite{Moi52} these two manifolds are diffeomorphic. In particular, $\mathfrak{R}(\mathcal{T}^3)=\textup{GL}_3(\Z)$.\\
\hf It is known (cf. \cite{Wall99} \S 15{\tiny A}) that $\mathcal{T}^n$ is homeomorphic to $\mathbb{T}^n$ and $|\textup{\Str}(\mathbb{T}^n)|=|\textup{\Str}(\mathcal{T}^n)|$ is finite when $n\geq 5$, whence $\textup{\Str}_0(\mathcal{T}^n)$ is finite as well. Using $h\mathfrak{R}(1,n)=\textup{GL}_n(\Z)$ (cf. \eqref{hRkn13}) and applying Corollary \ref{Rnk} to $\mathcal{T}^n$, we arrive at an analogue of Corollary \ref{Rknnor} for $\mathcal{T}^n$, i.e., $\mathfrak{R}(\mathcal{T}^n)$ contains an infinite normal subgroup of $\textup{GL}_n(\Z)$. Therefore, $\mathfrak{R}(\mathcal{T}^n)$ also contains a congruence subgroup. $\hfill\square$
\begin{rem}
It follows from the discussion in \S \ref{osph} that $\mathfrak{R}_{\textup{st}}(k,n)=\textup{GL}_n(\Z)$ for $k=1,3,7$ and $n\geq 1$. In Remark \ref{stndsp} we had seen that $\mathfrak{R}_{\textup{st}}(k,n)=h\mathfrak{R}_\textup{st}(k,n)$ when $k\neq 3,7$ is odd.
\end{rem}
\begin{rem}
When $n=2$, $\mathfrak{R}(k,2)$ certainly contains $E_2(2m!)$, i.e., the normal subgroup of $\textup{SL}_2(\Z)$ generated by the elementary matrices 
\bgd
E_{12}^{2m!}=\left(\begin{array}{cc}
1 & 2m!\\
0 & 1
\end{array}\right)\,\,\textit{and}\,\,E_{21}^{2m!}=\left(\begin{array}{cc}
1 & 0\\
2m! & 1
\end{array}\right),
\edd
where $m=\big|\textup{\Str}_0(\Sigma_1^k\times \Sigma_2^k)\big|$. However, $E_2(s)$ has infinite index in $\textup{SL}_2(\Z)$ whenever $s\geq 6$; cf. \cite{New67}. But
\bgd
\left(\begin{array}{cc}
1 & s\\
s & s^2+1
\end{array}\right)\in E_2(s)
\edd 
and is a hyperbolic matrix (i.e., its eigenvalues are not $1$ in absolute value) whenever $s\geq 2$.
\end{rem}
{\bf Proof of Lemma \ref{ssfin}.}\hf This result is a consequence of the Browder-Novikov-Sullivan-Wall surgery exact sequence\footnote{The $L$-groups are that of $L_\ast(\Z[\pi_1(M)])=L_\ast(\Z)$ as $M$ is simply connected. The map $\delta$ is actually induced by an action of $L_{nk+1}(\Z)$ on $\textup{\Str\,}(M)$ by connect sum with elements of $b\mathcal{P}_{nk+1}$.} of pointed sets \cite{Wall99}, i.e.,  
\bge\label{sur}
\xymatrix{
& & [M,\textup{G}]\ar[d]^-{\omega} & \\
L_{nk+1}(\Z)\ar[d]\ar[r]^-{\delta} & \textup{\Str\,}(M^{nk})\ar[r]^-{\eta} & [M,\textup{G}/\textup{O}]\ar[r]\ar[d]^-{\tau} & L_{nk}(\Z)\\
b\mathcal{P}_{nk+1}\ar[ru] & & [M,B\textup{O}]=\widetilde{\textup{KO}}(M) & }
\ede
where $b\mathcal{P}_{nk+1}$ is a subgroup of $\Theta_{nk}$. Here $\Theta_{nk}$ is the Kervaire-Milnor group of homotopy spheres of dimension $nk$ and $b\mathcal{P}_{nk+1}$ is the subgroup of exotic spheres that bound parallelizable manifolds. The main theorem of the Kervaire-Milnor paper \cite{KM63} states that $\Theta_{nk}$ is a finite group. Now the composition map $\tau\circ\eta$ in exact sequence \eqref{sur} can be described as follows. Let $f:N\to M$ represent an element $x\in\textup{\Str\,}(M)$, then $\tau(\eta(x))$ is represented by the stable vector bundle 
\bgd
\tau(\eta(x))=({\bf f}^\ast)^{-1}( f^\ast TM - TN)=TM-(f^{-1})^\ast(TN)\in \widetilde{\textup{KO}}(M).
\edd
In this formula ${\bf f}^\ast:\widetilde{\textup{KO}}(M)\to \widetilde{\textup{KO}}(N)$ is the functorially induced isomorphism and $TM,TN$ are the stable tangent bundles. Next note that $M=\Sigma_1^k\times\cdots\times \Sigma_n^k$ is stably parallelizable; so when $M=N$, we have
\bgd
\tau\big(\eta(x)\big)=0.
\edd
Consequently, 
\bgd
\eta\big(\textup{\Str}_0(M)\big)\subseteq \textup{image}\,\,\omega
\edd
where $\omega:[M,\textup{G}]\to [M,\textup{G}/\textup{O}]$ is induced by the natural map $\textup{G}\to \textup{G}/\textup{O}$. But $[M,\textup{G}]$ is a finite group as $M$ is a finite CW complex. As a consequence, $\textup{\Str}_0(M)$ is a finite set. $\hfill\square$
\begin{rem}\label{mainthm3}
One can improve Theorem \ref{mainthm} to include smooth manifolds $M$ which are homeomorphic to $n$-fold product of $k$-dimensional spheres, i.e., for such a manifold $M$, the subgroup $\mathfrak{R}(M)$ of $\textup{GL}_n(\Z)$ realizable by self-diffeomorphisms of $M$ contains a congruence subgroup. To prove this we need only show that $\textup{\Str}_0(M)$ is finite. Equivalently, it suffices to show that $\tau\big(\eta(\textup{\Str}_0(M))\big)$ has finite image. Since $M(k,n)$ is parallelisable, all higher Pontrjagin classes vanish. By Novikov's topological invariance of rational Pontrjagin classes, $p_i(M)\otimes\Q$ are zero for $i\geq 1$. But $\tau\big(\eta(\textup{\Str}_0(M))\big)$ is contained in the subset of stable bundles with vanishing (higher) rational Pontrjagin classes. Since the Pontrjagin character
\bgd
\textup{ch}:\widetilde{KO}(M)\lra \oplus_{i=0}^\infty H^{4i}(M;\Q)
\edd
is a rational isomorphism, it follows that elements of $\widetilde{KO}(M)$ with vanishing rational Pontrjagin classes are torsion elements. As $M$ is a compact manifold, $\widetilde{KO}(M)$ is a finitely generated abelian group and the torsion subgroup is finite. This implies that $\tau\big(\eta(\textup{\Str}_0(M))\big)$ is finite. 
\end{rem}


\bibliographystyle{siam}

\vspace*{0.4cm}

\hf {\small D}{\scriptsize EPARTMENT OF }{\small M}{\scriptsize ATHEMATICS, }{\small RKM V}{\scriptsize IVEKANANDA }{\small U}{\scriptsize NIVERSITY, }{\small H}{\scriptsize OWRAH, }{\small WB} {\footnotesize 711202, }{\small INDIA}\\
\hf{\it E-mail address} : \texttt{somnath@maths.rkmvu.ac.in}

\vspace*{0.4cm}

\hf {\small Y}{\scriptsize AU }{\small M}{\scriptsize ATHEMATICAL }{\small S}{\scriptsize CIENCES }{\small C}{\scriptsize ENTER, }{\small T}{\scriptsize SINGHUA }{\small U}{\scriptsize NIVERSITY, }{\small BEIJING, CHINA}\\
\hf{\it E-mail address} : \texttt{farrell@math.tsinghua.edu.cn}


\end{document}